\pdfoutput=1
\RequirePackage{ifpdf}
\ifpdf 
\documentclass[pdftex]{sigma}
\else
\documentclass{sigma}
\fi

\usepackage{bbm}
\usepackage{texdraw}
\usepackage[all]{xy}
\usepackage{pb-diagram}

\numberwithin{equation}{section}

\newtheorem{Theorem}{Theorem}[section]
\newtheorem*{Theorem*}{Theorem}
\newtheorem{Corollary}[Theorem]{Corollary}
\newtheorem{Lemma}[Theorem]{Lemma}
\newtheorem{Proposition}[Theorem]{Proposition}

\newtheorem*{Claim}{Claim}
{ \theoremstyle{definition}
 \newtheorem{Definition}[Theorem]{Definition}

 \newtheorem{Remark}[Theorem]{Remark}
 
}

\def\={ = }
\def\+{ + }

\def\be{\begin{equation}}
\def\ee{\end{equation}}

\begin{document}
\allowdisplaybreaks

\newcommand{\arXivNumber}{2401.05734}

\renewcommand{\PaperNumber}{027}

\FirstPageHeading

\ShortArticleName{The Morse--Smale Property of the Thurston Spine}

\ArticleName{The Morse--Smale Property of the Thurston Spine}

\Author{Ingrid IRMER}

\AuthorNameForHeading{I.~Irmer}

\Address{International Center for Mathematics, Department of Mathematics,\\ Southern University of Science and Technology, Shenzhen, P.R.~China}
\Email{\href{mailto:ingridmary@sustech.edu.cn}{ingridmary@sustech.edu.cn}}
\URLaddress{\url{https://www.sustech.edu.cn/en/faculties/ingrid-irmer.html}}

\ArticleDates{Received February 07, 2024, in final form April 10, 2025; Published online April 23, 2025}

\Abstract{The Thurston spine consists of the subset of Teichm\"uller space at which the set of shortest curves, the systoles, cuts the surface into polygons. The systole function is a topological Morse function on Teichm\"uller space. This paper studies the local properties of the Thurston spine, and the smooth pieces out of which it is constructed. Some of these local properties are shown to have global consequences, for example that the Thurston spine satisfies properties defined in terms of the systole function analogous to that of Morse--Smale complexes of (smooth) Morse functions on compact manifolds with boundary.\looseness=1}

\Keywords{Thurson spine; moduli space; handle decomposition; systole; topological Morse function; unstable manifold; mapping class group}

\Classification{57R99; 57M99}

\section{Introduction}
\label{sect:intro}

The systoles on a closed, connected, orientable hyperbolic surface $S_{g}$ of genus $g$ are the set of shortest geodesics. This is a finite nonempty set that varies as the hyperbolic structure varies in Teichm\"uller space $\mathcal{T}_{g}$. The Thurston spine $\mathcal{P}_{g}$ is a mapping class group-equivariant CW complex embedded in $\mathcal{T}_{g}$. Setwise it consists of the points of $\mathcal{T}_{g}$ where the systoles cut the corresponding surface into polygons. In \cite{Thurston}, a construction of a mapping class group-equivariant deformation retraction of $\mathcal{T}_{g}$ onto $\mathcal{P}_{g}$ is given. Although Morse theory was not explicitly used, the construction is analogous to the deformation retraction of a compact manifold with boundary onto the Morse--Smale complex defined by a Morse function.

A key role is played by the systole function $f_{\mathrm{sys}}\colon \mathcal{T}_{g}\rightarrow \mathbb{R}_{+}$ taking a point in $\mathcal{T}_{g}$ to the length of the systoles. The systole function is a piecewise-smooth topological Morse function, \cite{Akrout,SchmutzMorse}.

The purpose of this paper is to relate the structure of $\mathcal{P}_{g}$ to that of the equivariant handle decomposition of $\mathcal{T}_{g}$ coming from $f_{\mathrm{sys}}$. This is done by proving three results, outlined below.

To begin with, all critical points of $f_{\mathrm{sys}}$ are contained in $\mathcal{P}_{g}$~\cite{Thurston}. As $f_{\mathrm{sys}}$ is not smooth, the notion of gradient of $f_{\mathrm{sys}}$ at $x$ is replaced by a ``cone of increase'' consisting of the elements of~$T_{x}\mathcal{T}_{g}$ specifying directions in which $f_{\mathrm{sys}}$ is increasing. Similarly, as $\mathcal{P}_{g}$ is not a manifold, the notion of tangent space to $\mathcal{P}_{g}$ at a point $x$ is replaced by the ``tangent cone'' to $\mathcal{P}_{g}$ at $x$. The following theorem states informally that restricting $f_{\mathrm{sys}}$ to $\mathcal{P}_{g}$ does not create critical points.

\begin{Theorem}
Suppose a point $x\in \mathcal{P}_{g}$ has the property that the cone of increase of $f_{\mathrm{sys}}$ at $x$ is full. Then the intersection of this cone with the tangent cone of $\mathcal{P}_{g}$ is nonempty.
\end{Theorem}

In other words, at a point at which $f_{\mathrm{sys}}$ is increasing in some direction to first order, there exists a direction in the tangent cone to $\mathcal{P}_{g}$ in which the lengths of all the systoles are increasing.\looseness=1

The next two results ensure that the Thurston spine contains a choice of unstable manifolds\footnote{As the systole function is only a topological Morse function, the stable/unstable manifolds are not uniquely defined. This is one reason for wanting to work with the Thurston spine instead of the unstable manifolds. The notion of stable manifolds will be replaced by Schmutz's notion of sets of minima, as made rigorous in an upcoming paper~\cite{STduality}.} of the critical points of $f_{\mathrm{sys}}$.

\begin{Theorem}
Where a stratum of $\mathcal{P}_{g}$ is locally top-dimensional it is balanced. If $x\in \mathcal{P}_{g}$ is not balanced, there are points directly above $x$ in $\mathcal{P}_{g}$ that are balanced.
\end{Theorem}

Technical terms, such as stratum, tangent cone and cone of increase, are defined in Section \ref{sub:defns}. A balanced stratum has properties in common with an unstable manifold of the systole function, in the way the systole function increases towards such a stratum.

\begin{Proposition}
\label{thm3}
Fix a critical point $p\in \mathcal{P}_{g}$ of index $j$ with set of systoles $C$. A neighbourhood of $p$ in $\mathcal{P}_{g}$ contains a piecewise smooth cell of codimension $j$ on which $f_{\mathrm{sys}}$ is increasing to second order away from the critical point. The tangent space of this cell at $p$ has codimension $j$ and is given by $\{\nabla L(c)(p) \mid c\in C\}^{\perp}$. The subspace spanned by $\{\nabla L(c)(p) \mid c\in C\}$ consists of directions in which the systole function is decreasing away from $p$.
\end{Proposition}

The proofs given in this paper do not assume Thurston's equivariant deformation retraction of $\mathcal{T}_{g}$ onto $\mathcal{P}_{g}$, and are elementary in the sense that they only require Thurston's Lipschitz map construction from \cite{Thurston}, Wolpert's formula for Dehn twists, convexity of length functions, calculus and linear algebra. They will be used to prove geometric properties of gradients of length functions, to justify part of the construction of \cite{Thurston}, and in forthcoming papers to relate geometric, topological, analytic, algebraic and number theoretic properties of Teichm\"uller space and the mapping class group.

As $f_{\mathrm{sys}}$ is not smooth, there can potentially exist regular points at which there is no open cone of directions in which the systole function is increasing. These points are called \textit{boundary points} of the systole function. A geometric model for $\mathcal{P}_{g}$ around critical points and boundary points of $f_{\mathrm{sys}}$ is given in Section \ref{criticalstructure} and an explicit example of a critical point is worked through in Section \ref{subexample}. Whether or not boundary points actually exist in Teichm\"uller space of surfaces without boundary is not known. An example of a boundary point of $f_{\mathrm{sys}}$ for a surface with boundary is given in \cite[p.~434, point~(ii)]{SchmutzMorse}.

The Morse theory used in this paper is not the usual (smooth) Morse theory, but the Morse theory of so-called \textit{topological Morse functions} (defined in Section \ref{sub:defns}) developed in \cite{Morse}. It was shown in \cite{Akrout} that~$f_{\mathrm{sys}}$ is a topological Morse function. The critical points necessarily occur where~$f_{\mathrm{sys}}$ is not smooth. Studying the pieces of the piecewise smooth structure is a key ingredient of most of the proofs in this paper. These pieces, the strata, are contained in so-called ``loci''. These are solutions to the set of constraints that the lengths of curves in a given set are equal, or differ by a fixed amount encoded in a tuple $d$. Of particular interest are the loci described by sets of curves that do not contain any proper filling subsets; ``minimal filling sets''.\looseness=1

\begin{Corollary}
\label{thecor}
Let $C$ be a minimal filling set of curves. When the locus $E(C,d)$ is nonempty, it is a connected, embedded submanifold of $\mathcal{T}_{g}$. The restriction to $E(C,d)$ of the lengths of the curves in $C$ is strictly convex and has a unique minimum.
\end{Corollary}

Typically, locally top-dimensional strata of $\mathcal{P}_{g}$ are open sets of loci described by Corollary~\ref{thecor}, and important properties of such loci are shown to generalise to locally top-dimensional strata of $\mathcal{P}_{g}$.

The reference to Morse--Smale in the title comes from a CW-complex embedded in a manifold with boundary that is often referred to as the ``Morse--Smale complex'' in computer science literature. This complex is obtained as the image of a deformation retraction constructed using the flow of the gradient vector field of a Morse function. The cells of the complex consist of the unstable manifolds of the critical points. This CW-complex should not be confused with the chain complex of the same name used for computing homology of the manifold. An analogue of the Morse--Smale complex for topological Morse functions is described in~\cite{Morse}.\looseness=1

\textbf{Organisation of the paper.} In an effort to be self-contained, a fairly detailed Section~\ref{sub:defns} is given, surveying the basic properties of the objects studied in this paper and introducing the notation and assumptions that will be used. This includes a subsection on Schmutz's sets of minima; a construction that is perhaps not as well-known in the literature as it deserves to be. In order to understand the properties of the stratification of $\mathcal{T}_{g}$, the topological properties of the loci are studied in Section \ref{Minimalfilling}. An important concept appears to be that of a balanced locus, which is introduced in Section \ref{subbalanced}. Section \ref{mainproof} proves the two theorems stated at the beginning of the introduction, and similar properties are established for critical and boundary points of $f_{\mathrm{sys}}$ in Section \ref{criticalstructure}. The implications of the theorems in Sections \ref{mainproof} and \ref{criticalstructure} for the structure of a vector field studied by Thurston are discussed in Section \ref{Tvf}. In Section \ref{subexample}, Proposition \ref{thm3} is applied to an important example of a critical point due to Schmutz.

\section{Definitions and conventions}
\label{sub:defns}

The unique closed topological surface of genus $g$ for $g\geq 2$ will be denoted by $\mathcal{S}_{g}$. When $\mathcal{S}_{g}$ has the additional structure of a marked hyperbolic structure, i.e., it represents a point of Teichm\"uller space $\mathcal{T}_{g}$ of $\mathcal{S}_{g}$, it will be denoted by $S_{g}$.

\textit{Curves} on $\mathcal{S}_{g}$ are nontrivial isotopy classes of embeddings of $S^{1}$ into $S_{g}$, where $S^{1}$ is the unpointed, unoriented 1-sphere. A curve on $S_{g}$ is assumed to inherit a marking, and by abuse of notation, the word curve on $S_{g}$ will often be confused with the image of its unique geodesic representative. The term geodesic will always refer to a curve except in the case of a geodesic arc. The set of curves on $\mathcal{S}_{g}$ will be denoted by $\mathcal{C}_{g}$. A curve will usually be denoted by lowercase~$c$ with a subscript, while uppercase $C$ is reserved for finite sets of curves. A geodesic on $S_{g}$ of smallest length will be called a \textit{systole}. The \textit{systole function} $f_{\mathrm{sys}}$ is the piecewise smooth function~${\mathcal{T}_{g}\rightarrow \mathbb{R}_{+}}$ whose value at a point $x$ is given by the length of the systoles. Throughout the paper, it will always be assumed that the elements of a set $C$ of curves have pairwise geometric intersection number at most one, mirroring a property of sets of systoles.

A set of curves on $\mathcal{S}_{g}$ \textit{fills} the surface if the complement of a set of representatives in minimal position consists of a set of topological disks. The Thurston spine $\mathcal{P}_{g}$ is therefore the set of all points in $\mathcal{T}_{g}$ at which the systoles fill.

In this paper, a cell is homeomorphic to an open ball. The standard definition of a cell in a CW complex is closed; these will be understood to be the closures of the cells discussed in this paper. By abuse of notation, vertices, edges and faces of a cell will be used to refer to the vertices, edges and faces of the closure of the cell.

The symbol $\subset$ will be used to denote subset. For example, $\{1\}\subset \{1,2\}$ and $\{1,2\}\subset \{1,2\}$, whereas $\subsetneq$ is used to rule out the possibility that the sets are equal.

\textbf{Length functions.} Length functions generalise the notion of the length of a curve. A~curve~$c$ on $\mathcal{S}_{g}$ determines an analytic function $L(c)\colon\mathcal{T}_{g}\rightarrow \mathbb{R}_{+}$ whose value at a point $x$ is equal to the length of the geodesic representative of $c$ at $x$. For an ordered set of curves~${C=(c_{1}, \dots, c_{k})}$ and a~corresponding set of real, positive weights $A=(a_{1}, \dots, a_{k})$, the length function $L(A,C)\colon\allowbreak {\mathcal{T}_{g}\rightarrow \mathbb{R}_{+}}$ is the analytic function given by
\begin{equation*}
L(A,C)(x) = \sum_{j=1}^{k} a_{j}L(c_{j})(x).
\end{equation*}
When all the weights in $A$ are equal to 1, the map $L(A,C)$ from $\mathcal{T}_{g}$ to $\mathbb{R}^{C}$ will be written $L(C)$ or $L(c)$ when $C=\{c\}$. Length functions are so useful because they satisfy many convexity properties, as shown in \cite{Wolpert} for the Weil--Petersson metric, in \cite{Kerckhoff} for earthquake paths and \cite{Bestvina} for Fenchel--Nielsen coordinates.

The \textit{sublevel set} of $L(c)$ is the set of all points of $\mathcal{T}_{g}$ for which $L(c)$ is less than or equal the constant $l \in \mathbb{R}_{+}$, and will be denoted by $L(c)_{\leq l}$. The boundary of $L(c)_{\leq l}$ is then the level set~$L(c)^{-1}(l)$. Similarly, the superlevel set $\{x\in \mathcal{T}_{g} \mid L(c)(x)\geq l\}$ will be denoted by $L(c)_{\geq l}$.

\textbf{Local finiteness.} The lengths of elements of $\mathcal{C}$ is known to satisfy a property that will be referred to as local finiteness. It can be proven from the collar lemma that for any $x\in \mathcal{T}_{g}$, the set of curves of length within $\epsilon$ of $f_{\mathrm{sys}}(x)$ is finite. Consequently, there is a neighbourhood of~${x\in \mathcal{T}_{g}}$ in which the systoles are contained in the set of curves that are realised as systoles at~$x$.

\begin{Definition}[systole stratum]
A \textit{systole stratum} $\mathrm{Sys}(C)$ is the subset of $\mathcal{T}_{g}$ on which $C$ is the set of systoles.
\end{Definition}

Here the term stratification is used in a minimal sense, following Thurston, to mean a disjoint union of locally closed subsets, for which each point is on the boundary of only finitely many of these locally closed subsets (strata).

The existence of a triangulation of $\mathcal{T}_{g}$ compatible with the stratification follows from \cite{Lojasiewicz1964}. A triangulation of this type is assumed when referring to the boundary of a stratum, or the tangent vector or 1-sided limit of tangent vectors to strata of $\mathcal{P}_{g}$.

At first sight it might seem surprising that $\mathcal{P}_{g}$ is nonempty. This is a consequence of the following proposition, and of the theorem due to Bers' \cite{Bers1985} that $f_{\mathrm{sys}}$ is bounded from above.

\begin{Proposition}[{\cite[Proposition 0.1]{Thurston}, \cite[Lemma 4]{Bers}}]
\label{Thurstonprop}
Let $C$ be a set of curves on $\mathcal{S}_{g}$ that do not fill. Then at any point $x$ of $\mathcal{T}_{g}$, there is an open cone of derivations in $T_{x}\mathcal{T}_{g}$ whose evaluation on every element of $\{L(c) \mid c\in C\}$ is strictly positive at $x$.
\end{Proposition}

The following is a short proof of Proposition \ref{Thurstonprop} suggested by an anonymous referee.

\begin{proof}
Since $C$ does not fill, there exists a curve $c^{*}$ disjoint from all the curves in $C$. Dehn twisting around $c^{*}$ takes $x$ to a different point $x'$ without changing the lengths of curves in $C$. Since length functions are convex along Weil--Petersson geodesics~\cite{Wolpert}, the tangent vector $v$ at~$x$ to the Weil--Petersson geodesic from~$x$ to $x'$ is a direction in which the lengths of all the curves in $C$ are decreasing to first order. It follows that $v$ is contained in an open cone of vectors satisfying the proposition.
\end{proof}

\textbf{Twist flow.} In addition to a length function $L(c)$, a curve $c$ on $\mathcal{S}_{g}$ determines a flow $\mathbb{R}\times \mathcal{T}_{g}\rightarrow \mathcal{T}_{g}$. Choose any set of Fenchel--Nielsen coordinates on $\mathcal{T}_{g}$ defined using a set of curves containing~$c$. Then a point $(r,x)$ is mapped to the point in $\mathcal{T}_{g}$ with Fenchel--Nielsen coordinates equal to the coordinates at $x$, except for the twist parameter around $c$, which is incremented by~$r$. A~flowline will be called a \textit{twist path} and denoted by $\gamma_{c, x}$, where $x=\gamma_{c,x}(0)$. A twist path is a special case of an earthquake path, for which it is known that length functions are convex~\cite{Kerckhoff}.

This flow will be implicit in many of the arguments in this paper, for the simple reason that for any $c_{1}$ and $c_{2}$ disjoint, $c_{1}$ determines a flow that leaves the level sets of $L(c_{2})$ invariant.

\textbf{Topological Morse functions.} Suppose $M$ is an $n$-dimensional topological manifold. A~continuous function $f\colon M\rightarrow \mathbb{R}_{+}$ will be called a \textit{topological Morse function} if the points of~$M$ are either regular or critical points with respect to $f$. A point $x \in M$ is regular if there is an open neighbourhood $U$ containing $x$, for which $U$ admits a homeomorphic parametrisation by $n$ parameters, one of them being $f$. If the point $p$ is a critical point, there exists a $j\in \mathbb{Z}$, $0\leq j\leq n$, called the \textit{index} of $p$, and a homeomorphic parametrisation of $U$ by parameters~${\{x_{1}, \dots, x_{n}\}}$ such that for every point in~$U$
\[
f(x) - f(p) = \sum_{i=1}^{i=n-j}x_{i}^{2} - \sum_{i=n-j+1}^{i=n}x_{i}^{2}.
\]

A reference for topological Morse functions is \cite{Morse}.

Whenever a metric is needed in this paper, the Weil--Petersson metric will be assumed, and denoted by $\mathbf{g}$. None of the important concepts are dependent on a specific choice of mapping class group-equivariant metric; in most instances, any other choice would work just as well.

The orthogonal complement of a set of gradients of length functions, $\{\nabla L(c) \mid c\in C\}$ will be denoted by $\{\nabla L(c) \mid c\in C\}^{\perp}$. Although a metric is needed to define gradients and orthogonal complements, $\{\nabla L(c) \mid c\in C\}^{\perp}$ is independent of a choice of metric. This is because it can be defined as a tangent space to the intersection of level sets of length functions.

\textbf{Tangent cones and cones of increase.} There are many places throughout this paper at which objects similar to polyhedra are used. For a polyhedron $P$ in $\mathbb{R}^{d}$, the \textit{tangent cone} at a~point $x\in \partial P$ is the set of $v\in T_{x}\mathbb{R}^{d}$ such that $v=\dot{\gamma}(0)$ for a smooth oriented path $\gamma(t)$ with~${\gamma(0)=x}$ and $\gamma(\epsilon)$ in $P$ for sufficiently small $\epsilon>0$. For background and definitions relating to polyhedra, see \cite{BVcone} or \cite{Tcone}.

As $\mathcal{P}_{g}$ can be triangulated, the tangent cone of $\mathcal{P}_{g}$ at a point $x$ can be defined analogously. When $x$ is on the boundary of more than one simplex of $\mathcal{P}_{g}$, the tangent cone is the union of the tangent cones of the simplices with $x$ on the boundary. Tangent cones of strata are defined analogously. For a triangulation compatible with the stratification, the tangent cone to $\mathrm{Sys}(C)$ at $x$ is the union of the tangent cones of the simplices with interior contained in $\mathrm{Sys}(C)$ with $x$ on the boundary.

Since $f_{\mathrm{sys}}$ is only piecewise smooth, it does not have a gradient at every point. The level sets of $f_{\mathrm{sys}}$ are a locally finite intersection of level sets of curves in $\mathcal{T}_{g}$. On a neighbourhood of $x\in \mathcal{T}_{g}$ at which $C$ is the set of systoles, by local finiteness the level set of $f_{\mathrm{sys}}$ passing through $x$ lies on the boundary of the intersection $I(C,x)$ defined as follows
$
I(C,x):=\bigcap_{c\in C}L(c)_{\geq f_{\mathrm{sys}}(x)}$.
Each of these sublevel sets is linear with respect to some set of Fenchel--Nielsen coordinates, but as the systoles can intersect, $I(C,x)$ is not always an intersection of half-spaces linear with respect to the same set of coordinates. Nevertheless, it is possible to define a tangent cone in the same way as for a polyhedron, and the \textit{cone of increase} of $f_{\mathrm{sys}}$ at a point $x$ is defined to be the tangent cone of $I(C,x)$.

The cone of increase of $f_{\mathrm{sys}}$ will be called \textit{full} if it has dimension equal to that of $\mathcal{T}_{g}$, namely~${6g-6}$. Since $f_{\mathrm{sys}}$ is only a topological, but not smooth, Morse function, it is possible that there exist regular points at which the cone of increase of $f_{\mathrm{sys}}$ is not full. If $x$ is one such point, then $f_{\mathrm{sys}}$ can increase away from $x$ to at most second order. Following Schmutz, these points will be called \textit{boundary points of }$f_{\mathrm{sys}}$. It was shown in \cite{SchmutzMorse} that the set of critical and boundary points of $f_{\mathrm{sys}}$ are a discrete subset of $\mathcal{P}_{g}$. An alternative definition will be given in Section \ref{SchmutzSection}, and these points will be studied in Section \ref{criticalstructure}.

Wherever it is nonempty, $\partial I(C,x)$ is a piecewise smooth, codimension 1 submanifold of~$\mathcal{T}_{g}$, and as such, it makes sense to speak of $y\neq x$ on $\partial I(C,x)$. A point $y$ on $\partial I(C,x)$ will be called a \textit{vertex} of $I(C,x)$ if the tangent cone to $I(C,x)$ at $y$ contains no vector subspaces of~$T_{y}\mathcal{T}_{g}$. Otherwise, $y$ will be said to be on a \textit{face}, where the face corresponds to a maximal vector subspace of the tangent cone of $I(C,x)$ at $y$, and two such subspaces at points $y_{1}$ and $y_{2}$ determine the same face if $y_{1}$ and $y_{2}$ are connected by a path in $\partial I(C,x)$ along which the vector subspace has fixed dimension and varies smoothly. When the face has codimension 1 in $\mathcal{T}_{g}$, it will be called a \textit{facet}.

When $I(C,x)$ is replaced by the intersection
$
 D(C,t):=\bigcap_{c\in C}L(c)_{\leq t}
$
 tangent cones, facets, faces and vertices are defined analogously.

The fact that $f_{\mathrm{sys}}$ is a topological Morse function~\cite{Akrout} implies that the set of unit tangent vectors in the cone of increase at a regular point of $f_{\mathrm{sys}}$ is a topological ball~\cite{Morse}. Removing the zero vector from the cone of increase does not give a disconnected set.

\textbf{Eutactic.} In the literature, a point $x\in \mathcal{T}_{g}$ in $\mathrm{Sys}(C)$ is defined to be \textit{eutactic} if for every derivation $v$ in $T_{x}\mathcal{T}_{g}$, either the evaluations $\{vL(c)(x) \mid c \in C\}$ are all zero, or there is at least one $c\in C$ for which $vL(c)(x)$ is strictly positive and another for which it is strictly negative. In this paper the term eutactic will be generalised. A set of curves $C$ (not necessarily systoles) will be called eutactic at a point $x$ in $\mathcal{T}_{g}$ if for every derivation $v$ in $T_{x}\mathcal{T}_{g}$, either the evaluations~${\{vL(c)(x) \mid c \in C\}}$ are all zero, or there is at least one $c\in C$ for which $vL(c)(x)$ is strictly positive and another for which it is strictly negative. Similarly, for a fixed subspace $V\in T_{x}\mathcal{T}_{g}$, the set of curves $C$ will be called $V$-eutactic at $x$ if for every derivation $v$ in $V$, either the evaluations~${\{vL(c)(x) \mid c \in C\}}$ are all zero, or there is at least one $c\in C$ for which $vL(c)(x)$ is strictly positive and another for which it is strictly negative.

\begin{figure}
\centering
\includegraphics[width=0.94\textwidth]{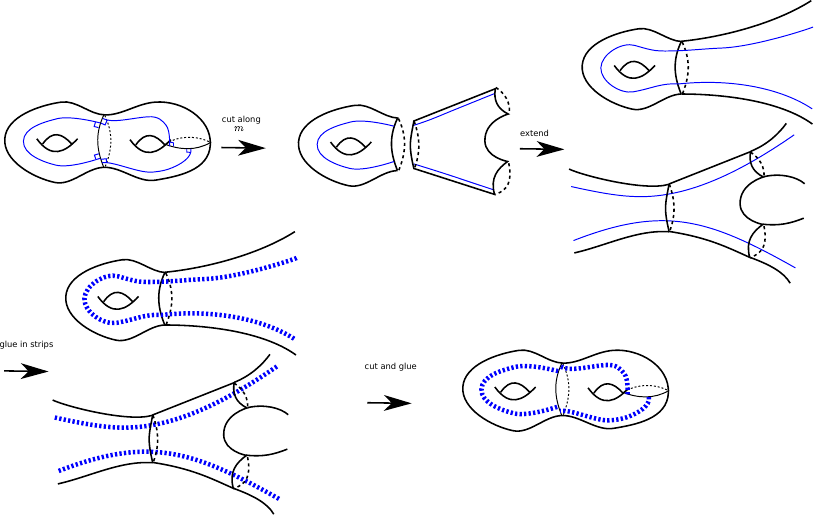}
\caption{Thurston's Lipschitz map construction, where $m$ is shown in black and $A$ in blue in the upper leftmost figure.}
\label{Lipschitzfig}
\end{figure}

\textbf{Thurston's Lipschitz map construction.} A construction of a map $\mathcal{T}_{g}\rightarrow \mathcal{T}_{g}$ will now be surveyed. This is explained in detail in \cite{PT}. The input of this construction is $(m,A)$, where $m$ is a geodesic multicurve on $S_{g}$ and $A$ is a set of pairwise disjoint geodesic arcs with endpoints on~$m$, meeting $m$ at right angles, and with the property that $A$ has nonzero geometric intersection number with every boundary component of a tubular neighbourhood of $m$. The output of the construction is a family of maps $\mathcal{T}_{g}\rightarrow \mathcal{T}_{g}$ which are used to control the change in the lengths of some curves relative to others. This construction was used in \cite{Thurston} to prove Proposition \ref{Thurstonprop}, with $m$ isotopic to the boundary of the subsurface filled by $C$, and $A$ a set of arcs that intersect every curve in $C$.

Start with a surface $S_{g}$ representing a point of $\mathcal{T}_{g}$. Cut $S_{g}$ along the geodesic representative of~$m$ to get a union of hyperbolic surfaces with boundary. On every boundary component of each of these surfaces, attach a hyperbolic annulus, to obtain a surface without boundary, as shown in Figure~\ref{Lipschitzfig}. Each of the geodesic representatives of arcs in $A$ can be extended to a bi-infinite geodesic in these extended surfaces. Cut the extended surfaces along each bi-infinite geodesic, and glue in a strip of the hyperbolic plane bounded by two infinite geodesics. The widths of the geodesic strips is chosen in such a way that for every curve $c$ in $m$, the curves $c_{L}$ and~$c_{R}$ in the extended surfaces corresponding to $c$ have the same length after the strips have been glued in. For each curve $c$ in $m$, cut the extended surfaces along the geodesic representatives of~$c_{L}$ and~$c_{R}$, and glue the surface together along the resulting boundary components. This gives a~new hyperbolic structure. Clearly there are many choices involved in this construction.

In this paper, as in \cite{Thurston}, limits are taken as the widths of the strips are allowed to approach zero. The resulting 1-parameter families of maps do not change the length of any curve disjoint from the geodesic representatives of the elements of $m$ and $A$. If the widths are chosen to increase to first order with the parameter, the tangent vectors to these 1-parameter families lengthen to first order every curve in $m$ and every curve with geodesic representative disjoint from $m$ that intersects one or more of the arcs in $A$.

\subsection{Sets of minima}
\label{SchmutzSection}

As an attempt at parameterising $\mathcal{T}_{g}$ without reference to punctures or marked points, Schmutz introduced the concept of cell decompositions of Teichm\"uller space parametrised by length functions. Schmutz's construction will be surveyed here as it provides insight into the differential topology of $\mathcal{T}_{g}$.

A length function $L(A,C)$ with fixed \smash{$A\in \mathbb{R}_{+}^{|C|}$} is convex with respect to the Weil--Petersson metric \cite{Wolpert}. A necessary and sufficient condition for $x$ to represent the minimum of $L(A,C)$ in~$\mathcal{T}_{g}$ is therefore that the gradient of $L(A,C)$ is zero at $x$. When $C$ is eutactic at $x$, there exists an $A$ for which the unique minimum of $L(A,C)$ occurs at $x$. If $C$ is not a filling set, Proposition~\ref{Thurstonprop} ensures that $L(A,C)$ has no minimum in the interior of $\mathcal{T}_{g}$.

\begin{Definition}[$\mathrm{Min}(C)$]
The subset $\mathrm{Min}(C)$ of $\mathcal{T}_{g}$ is the set of all points at which $L(A,C)$ has a minimum for some $A$ with strictly positive entries. Alternatively, $\mathrm{Min}(C)$ is the set of all~${x\in \mathcal{T}_{g}}$ at which $C$ is eutactic.
\end{Definition}

In Section \ref{Minimalfilling}, minima of sets of curves will be defined. This should not be confused with minima of length functions in the definition of $\mathrm{Min}(C)$.

There is a notion of closure of $\mathrm{Min}(C)$ that will be denoted by $\overline{\mathrm{Min}(C)}$.

\begin{Definition}[{modified version of \cite[Proposition 2]{SchmutzVoronoi}}]\label{defnmod}
Let $C$ be a finite set of closed geodesics that fills $S_{g}$. A point $x$ of $\mathcal{T}_{g}$ is in $\overline{\mathrm{Min}(C)}$ if and only if there does not exist a derivation in $T_{x}\mathcal{T}_{g}$ whose evaluation on each length function of a curve in $C$ is strictly positive. Then $\partial \mathrm{Min}(C)$ is defined to be the set of points in $\overline{\mathrm{Min}(C)}\setminus \mathrm{Min}(C)$
\end{Definition}

\begin{Lemma}[{consequence of \cite[Lemma 14 and Theorem 15]{SchmutzMorse}}]
\label{lemma14}
Any point $x\in\mathcal{T}_{g}$ in $\partial\mathrm{Min}(C)$ is in $\mathrm{Min}(C_{i})$, for some filling $C_{i}\subsetneq C$. If $C_{i}$ is the largest subset of $C$ that is eutactic at $x$, there exists a derivation $v$ in $T_{x}\mathcal{T}_{g}$ for which $vL(c)(x)=0$ for all $c\in C_{i}$ and $vL(c)(x)<0$ for all $c\in C\setminus C_{i}$.
\end{Lemma}

Lemma 1 of \cite{SchmutzMorse} showed that when $C$ fills, any choice of $A$ with strictly positive entries determines a length function $L(A,C)$ with a minimum realised somewhere in $\mathcal{T}_{g}$. It therefore follows from Proposition~\ref{Thurstonprop} that $\mathrm{Min}(C)$ is nonempty if and only if $C$ is a filling set of curves. There exists therefore a surjective map \smash{$\phi_{C}\colon\mathbb{R}^{|C|}_{+}\twoheadrightarrow \mathrm{Min}(C)$}, where $\phi_{C}(A)$ is the point in $\mathrm{Min}(C)$ at which $L(A,C)$ has its minimum. Since $\phi_{C}$ is not injective, the entries of the $k$-tuple $A$ are parameters, not coordinates. For $C=\{c_{1}, \dots, c_{k}\}$, let $F(C)\colon\mathcal{T}_{g}\rightarrow \mathbb{R}^{k}_{+}$ be the differentiable function $x\mapsto \bigl(L(c_{1}), \dots, L(c_{k})\bigr)$. When the rank of $F(C)$ is constant on $\mathrm{Min}(C)$, $\mathrm{Min}(C)$ is a~cell, otherwise it is ``pinched'' at places where the rank drops.

For a filling set $C$, some points on the boundary of $\mathbb{R}^{k}_{+}$ also determine length functions with minima in $\mathcal{T}_{g}$; Schmutz called these admissible boundary points. The union of $\mathbb{R}_{+}^{k}$ with the admissible boundary points will be denoted by $A_{C}$. In \cite[Lemma 2]{SchmutzMorse}, it was shown that there exists an open set $B_{C}$ in $\mathbb{R}^{k}$ containing $A_{C}$, with the property that every tuple $(b_{1}, \dots, b_{k})$ in~$B_{C}$ gives a function $b_{1}L(c_{1})+\dots+b_{k}L(c_{k})$ with a unique minimum in $\mathcal{T}_{g}$. The image under $\phi_{C}$ of~${B_{C}\setminus A_{C}}$ will be denoted by $\mathrm{Conv}(C)$, where ``$\mathrm{Conv}$'' stands for convex.

\textbf{Boundary points of the systole function.} It was shown in \cite{Akrout} that a critical point of $f_{\mathrm{sys}}$ occurs exactly where, for some set $C$, $\mathrm{Sys}(C)$ intersects $\mathrm{Min}(C)$. Schmutz defined a boundary point of $f_{\mathrm{sys}}$ to be a point at which, for some $C$, $\mathrm{Sys}(C)$ intersects $\partial\mathrm{Min}(C)$. By Lemma \ref{lemma14}, this is a point at which $\mathrm{Sys}(C)$ intersects $\mathrm{Min}(C')$, for $C'\subsetneq C$. There are reasons for calling such points boundary points of $f_{\mathrm{sys}}$, but boundary points of $f_{\mathrm{sys}}$ should not be confused with the usual topological notion of boundary. The structure of critical and boundary points of $f_{\mathrm{sys}}$ will be discussed in Section \ref{criticalstructure}. By \cite[Corollary 22, Theorem 23]{SchmutzMorse}, critical points and boundary points of $f_{\mathrm{sys}}$ are isolated, so there are only finitely many of them modulo the action of the mapping class group.

\section{Minimal filling sets of curves}
\label{Minimalfilling}

Suppose $C$ is a filling set of curves on $S_{g}$, any two of which intersect in at most one point. The set $C$ will be called a \textit{minimal filling set of curves} if no proper subset of $C$ fills. This section begins by proving some key lemmas about minimal filling sets of curves, and ends by discussing a property that generalises many of the important properties of minimal filling sets of curves.

A basic property of minimal filling sets of curves is the following, which will be assumed from now on without comment.

\begin{Lemma}
Suppose $C$ is a minimal filling set of curves. For any curve $c\in C$, there is a~curve~$c^{*}$ that intersects $c$ but does not intersect any of the other curves in $C$.
\end{Lemma}
\begin{proof}
Suppose every curve on $\mathcal{S}_{g}$ that intersects $c\in C$ also intersects a curve in $C\setminus \{c\}$. Then~${C\setminus \{c\}}$ fills, contradicting the assumption that $C$ is a minimal filling set.
\end{proof}

All known examples of locally top-dimensional strata of $\mathcal{P}_{g}$ have systoles that are minimal filling sets. The difficulty in proving this is always the case is due to a very specific configuration of linearly dependent length functions discussed in Remark \ref{thedifficulty}. It will be shown that loci of minimal filling sets of curves have beautiful properties that generalise to locally top-dimensional strata of $\mathcal{P}_{g}$.

\subsection{Key lemmas}
This subsection proves some key lemmas about minimal filling sets of curves and their loci, which will now be defined.

\begin{Definition}[locus $E(C)$ labelled by $C$]
The locus of points in $\mathcal{T}_{g}$ at which the lengths of all curves in $C$ are equal will be denoted by $E(C)$, and called the locus labelled by $C$. Alternatively,
$E(C)=\{x\in \mathcal{T}_{g} \mid L(c)(x)=L(c')(x) \, \forall c,c'\in C\}$.
\end{Definition}
Note that $\mathrm{Sys}(C)\subset E(C)$.

The set $E(C)$ is called a locus to emphasise the similarity with the locus of points equidistant from a finite set of points $\{p_{1}, \dots, p_{k}\}$ in $\mathbb{R}^{n}$. The analogy is drawn by replacing the set of points in $\mathbb{R}^{n}$ by points in the completion of $\mathcal{T}_{g}$ with respect to the Weil--Petersson metric, where $p_{i}$ for~${i=1, \dots, k}$ is a point at which the curve $c_{i}\in C$ has been pinched. Note that by convexity, $L(c_{i})$ increases monotonically along a Weil--Petersson geodesic emerging from $p_i$.

Recall that $\mathcal{C}$ is the set of curves on $\mathcal{S}_{g}$. Denote by $F(\mathcal{C})$ the map from $\mathcal{T}_{g}$ to $\mathbb{R}^{\mathcal{C}}$ taking $x$ to the collection of lengths $\{L(c)(x) \mid c\in \mathcal{C}\}$, and by $\overline{F}(\mathcal{C})\colon\mathcal{T}_{g}\rightarrow \mathbb{R}^{\mathcal{C}}/\mathbb{R}$ the composition of this map with the projection from $\mathbb{R}^{\mathcal{C}}$ (arbitrary maps from $\mathcal{C}$ to $\mathbb{R}$) to its quotient by $\mathbb{R}$ (the constant maps). Fix $C$ and $d\in \mathbb{R}^{\mathcal{C}}/\mathbb{R}$. Then the locus $E(C,d)$ is defined to be the set of points in $\mathcal{T}_{g}$ for which the lengths of the curves in $C$ realise the specified differences in $d|_{C}$, in other words
\begin{equation}
\label{withd}
E(C,d)=\{x\in \mathcal{T}_{g} \mid L(c)(x)+d|_{c}=L(c')(x)+d|_{c'} \, \forall c,c'\in C\}.
\end{equation}

When every element of $d|_{C}$ is zero, $d$ will be omitted.

The next lemma will be used in proving results about loci of minimal filling sets of curves.

\begin{Lemma}\label{LipschitzMap}
For a filling set $C$ of curves, suppose $C\setminus \{c\}$ is not filling for some $c\in C$. Then for any $x\in \mathcal{T}_{g}$, there exists a derivation $v\in T_{x}\mathcal{T}_{g}$ such that $vL(c)(x)<0$ and $vL(c')(x)>0$ for all $c'\in C\setminus \{c\}$.
\end{Lemma}

\begin{Definition}[$C$-generic]
Suppose $C$ is a (not necessarily proper) subset of a minimal filling set of curves. A point $x\in \mathcal{T}_{g}$ is $C$-generic if for every $c\in C$, the sum of the cosines of the angles of intersection of $c$ with $c^{*}$ is nonzero, for some choice of $c^{*}$. Alternatively, the length function~$L(c)$ is strictly convex along the twist path $\gamma_{c^{*},x}$ and hence has a unique minimum. The point~$x$ is $C$-generic if for every $c\in C$, $x$ is not the unique minimum of $L(c)$ on $\gamma_{c^{*},x}$.
\end{Definition}

\begin{proof}
Wolpert's twist formula \cite{Wolperttwist} states that the rate of change of the length of $L(c)$ with respect to the twist parameter around $c^{*}$ is equal to the sum of the cosines of the angles of intersection of $c$ with $c^{*}$. For points of $\mathcal{T}_{g}$ that are $\{c\}$-generic, the lemma follows from Wolpert's twist formula and Proposition \ref{Thurstonprop}. For points that are not $\{c\}$-generic, the proof uses Thurston's Lipschitz map construction from Section \ref{sub:defns}. The geodesic multicurve $m$ is given by $c^{*}$, where, as usual, $c^{*}$ is a curve intersecting $c$ and disjoint from any curve in $C\setminus\{c\}$. It will now be explained how to obtain the set of arcs $A$.

If $c$ intersects $c^{*}$ just once, and the cosine of the Weil--Petersson angle of intersection is zero, then $c\cap (S_{g}\setminus c^{*})$ is an arc with endpoints on $c^{*}$ making an angle $\pi/2$ with $c^{*}$. Choose this arc to be an element of $A$ and add more arcs to the set, subject to the constraints in the definition of~$A$, until each curve in the set $C\setminus \{c\}$ intersects at least one arc in $A$. As explained below, this choice of $m$ and $A$ then determines a 1-parameter family of maps, with the property that the tangent to this 1-parameter family gives a derivation that is strictly negative when evaluated on~$L(c)$ and strictly positive when evaluated on the length of any curve in $C\setminus\{c\}$.

When $c$ intersects $c^{*}$ more than once, it is necessary to keep in mind that the intersections between the arcs in $A$ and the curve $c$ that govern whether or not a given Lipschitz map increases or decreases the length of $c$ are not the usual topological intersections, but intersections between the geodesic representative of the curve $c$, and the geodesic representatives of the homotopy classes of arcs relative to the geodesic representative of $c^{*}$, i.e., geodesic arcs orthogonal at their endpoints to the geodesic $c^{*}$.

Choose an orientation for $c^{*}$. This makes it possible to define a ``right side'' $\mathcal{N}_{R}$ of $c^{*}$ in a~tubular neighbourhood $\mathcal{N}$ of $c^{*}$ and a ``left side'' $\mathcal{N}_{L}$ of $c^{*}$ in $\mathcal{N}$. Moving around $c^{*}$ and measuring angles to the right of $c^{*}$, there must be an arc $e_{1}$ of $c^{*}\cap (S_{g}\setminus \{c\})$ at whose endpoints the signs of the cosines of the angles of intersection either go from $+$ to $-$, $0$ to $-$, $+$ to~$0$ or~$0$ to~$0$. Here the angles to the right of $c^{*}$ at the endpoints of the arc $e_{1}$ are taken to be the angles at the vertices of the connected component of $S_{g}\setminus \{c, c^{*}\}$ to the right of $e_{1}$. In Figure~\ref{ccstar}\,(a), these are the angles~${\pi-\theta_{1}}$ and $\theta_{2}$. Similarly, moving around $c^{*}$ and measuring angles to the left of~$c^{*}$, there must be an arc $e_{2}$ at which the signs of the cosines of the angles of intersection either go from $+$ to $-$, $0$ to $-$, $+$ to $0$ or $0$ to~$0$. The arc $e_{1}$ joints two adjacent vertices of a connected component of $S_{g}\setminus \{c, c^{*}\}$ to the right of $e_{1}$, for which the inner angles at the two vertices are both at most $\frac{\pi}{2}$; see Figure~\ref{ccstar}\,(a). The arc $e_{2}$ does the same for a connected component of~$S_{g}\setminus \{c, c^{*}\}$ to the left of~$e_{2}$. These arcs will be used to construct arcs on both sides of the geodesic~$c^{*}$, called central arcs, whose geodesic representatives are disjoint from the geodesic~$c$.\looseness=1

\begin{figure}[ht]
\centering
\includegraphics[width=0.5\textwidth]{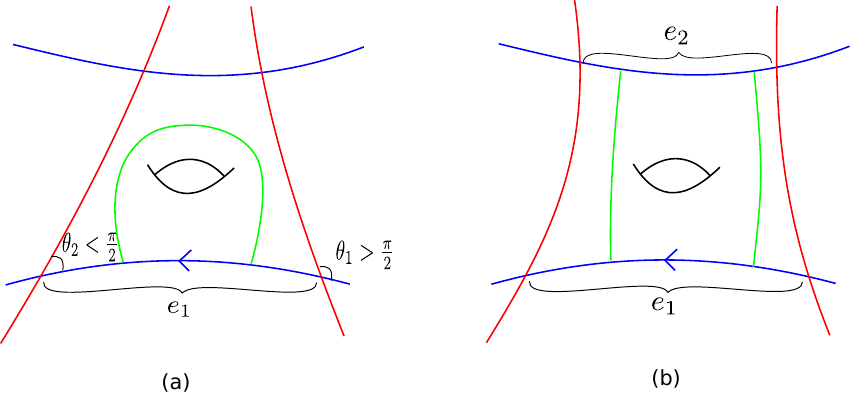}
\caption{In this figure, the blue arcs lie along the geodesic $c^{*}$ and the red arcs lie along the geodesic $c$ in the hyperbolic surface $S_{g}$. In (a), the green arc is a central arc disjoint from $c$ that only intersects $\mathcal{N}$ in $\mathcal{N}_{R}$, i.e., it is to the right of $c^{*}$. In (b), when $e_{1}$ and $e_{2}$ are contained in the same cycle, the central arcs are constructed as shown in green. These central arcs intersect both $\mathcal{N}_{R}$ and $\mathcal{N}_{L}$.}
\label{ccstar}
\end{figure}

A set of cycles $A_{i}:=\{a_{1}', \dots, a_{m}'\}$ is obtained, consisting of the connected components of the boundary of a normal neighbourhood of the union of the geodesics $c$ and $c^{*}$. By choosing~$c^{*}$ to have the smallest possible geometric intersection number with $c$, it is possible to assume without loss of generality that each element of $A_{i}$ is homotopically nontrivial. This is because the curves in the set $C$ intersect pairwise at most once, and $c$ is the only curve in $C$ to intersect~$c^{*}$. A curve in $C\setminus \{c\}$ that intersects a homotopically trivial element $a_{k}'$ of $A_{i}$ must intersect $c$ at least twice, which is not possible. When $a_{k}'$ is homotopically trivial, the curve~$c^{*}$ can then be replaced by a~homotopically nontrivial connected component of the boundary of a normal neighbourhood of~${c^{*}\cup a_{k}'}$ that intersects $c$.

Each cycle in $A_{i}$ has edges alternating between edges lying along $c$ and $c^{*}$. Suppose $e_{1}$ and~$e_{2}$ are not contained in the same cycle. From each cycle $a_{l}'$, remove a single edge lying along $c^{*}$ to obtain an arc $a_{l}$ with endpoints on $c^{*}$. Whenever possible, the arc removed should be one of the arcs $e_{1}$ or $e_{2}$. This is shown in Figure~\ref{ccstar}\,(a). An arc $a_{l}$ obtained by deleting~$e_{1}$ or $e_{2}$ will be called \textit{central}. When $e_{1}$ and $e_{2}$ are contained in the same cycle, a pair of central arcs are constructed as shown in Figure~\ref{ccstar}\,(b). A central arc has the property that the shortest representative of its homotopy class with endpoints on $c^{*}$ is disjoint from the geodesic representative of~$c$.

Choose a set $A$ of pairwise disjoint arcs as above, containing $\{a_{1}, \dots, a_{m}\}$.

If $A$ only contained central arcs, it follows that the resulting Lipschitz map construction gives a direction in which $L(c)$ is strictly decreasing and the lengths of curves in $C\setminus \{c\}$ are increasing. That the lengths of curves in $C\setminus \{c\}$ are increasing in this direction is a property of Thurston's Lipschitz map construction. Trigonometry implies that $L(c)$ decreases in this direction, as follows: Gluing in strips corresponding to the central arcs, as shown in Figure~\ref{newfig}, takes the geodesic~$c^{*}$ to a piecewise geodesic, with angles less than $\pi$ on the side of~$c^{*}$ that is not the flaring end. The geodesic representative of the piecewise geodesic is obtained by shifting the piecewise geodesic to the side with angles less than $\pi$, thereby decreasing~$L(c)$.

\begin{figure}[ht]
\centering
\includegraphics[width=0.5\textwidth]{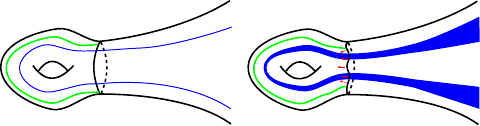}
\caption{The geodesic $c^{*}$, the central arc $a_{1}$ extended out into the flaring end (blue) and an arc of $c\cap (S_{g}\setminus c^{*})$ (green) are shown on the left. Gluing in a strip as shown turns $c^{*}$ into a piecewise smooth curve, which is shortened by shifting in the direction of the red arrows.}
\label{newfig}
\end{figure}

A Lipschitz map constructed from $m$ and $A$, for which the widths of the strips corresponding to the central arcs are sufficiently large relative to the widths of any strips corresponding to arcs that are not central, gives a 1-parameter family of maps whose tangent vector at $x$ corresponds to a derivation that is strictly negative when evaluated on $L(c)$ and strictly positive when evaluated on the length of any curve in $C\setminus\{c\}$. It is possible to simultaneously satisfy the constraints in Thurston's Lipschitz map construction and choose the strips corresponding to the central arcs to be wider than the other strips, because it was shown that there are at least~2 central arcs; one on either side of~$c^{*}$, or both intersecting both connected components of~$\mathcal{N}\setminus c^{*}$.
\end{proof}

Given a set $C$ of curves, and a point $x\in \mathcal{T}_{g}$, denote by $d(x)$ the value of $d$ for which~$E(C,d(x))$ passes through $x$. Choose $c\in C$, then a point in $E(C,d(x))$ at which $L(c)|_{E(C,d(x))}$ has a~minimum is a point at which the length of every curve in $C$, when restricted to $E(C,d(x))$, has a~minimum. A point $x\in \mathcal{T}_{g}$ will be called a \textit{minimum} of $C$ if $L(c)|_{E(C,d(x))}$ has a minimum at~$x$ for every $c\in C$.

\begin{Lemma}
\label{polyhedronargument}
Let $C$ be a minimal filling set of curves. Then the gradients of the lengths of the curves in $C$ are linearly independent everywhere in $\mathcal{T}_{g}$ away from minima of $C$.
\end{Lemma}
\begin{Remark}
Note that this lemma also holds for subsets of minimal filling sets of curves.
\end{Remark}

\begin{proof}
As in the proof of Lemma \ref{LipschitzMap}, there are two cases. When $x$ is $C$-generic, the lemma follows from Wolpert's twist formula. To prove the lemma for $C$-nongeneric points that are not minima of $C$, it will be shown that $C$-generic points satisfy a geometric property that implies linear independence of the gradients. It will then be shown that away from minima, Lemma \ref{LipschitzMap} implies this property can only break down on open sets. As the set of $C$-nongeneric points are not open, the lemma then follows.

Let $L(C, t, d):=\{x\in \mathcal{T}_{g} \mid L(c)(x)=t+d_{c} \, \forall c\in C\}$ be a level set of $E(C,d)$ from~\eqref{withd} on which the length of each $c\in C$ is given by $t+d|_{c}$. Then $D(C, t, d)$ is defined to be the intersection of the sublevel sets $\{L(c)_{\leq t+d|_{c}} \mid c\in C\}$. Recall the definition of faces, facets and tangent cones from Section \ref{sub:defns}. Since each sublevel set is strictly convex with respect to, for example, the Weil--Petersson metric~\cite{Wolpert}, the intersection $D(C, t, d)$ is convex, in particular connected. Each facet of $D(C,t,d)$ is contained in a level set $L(c)^{-1}(t+d|_{c})$ for some $c\in C$. Lower-dimensional faces of $D(C, t, d)$ are intersections of level sets.

The assumption that $E(C)$ is nonempty guarantees that for sufficiently large $t$ the intersection~${L(C,t)=\{x\in \mathcal{T}_{g} \mid L(c)(x)=t \, \forall c\in C\}}$ will be nonempty, and is then the lowest-dimensional face of $D(C,t)$.

Let $\{x_{1}, \dots, x_{n}\}$ be a standard set of coordinates on $\mathbb{R}^{n}$, and denote by
\[
\mathfrak{H}^{n}=\{(x_{1}, \dots, x_{n})\in \mathbb{R}^{n} \mid x_{1}x_{2}\cdots x_{n}=0\}\subset\mathbb{R}^{n}
\]
a union of hyperplanes passing through the origin of $\mathbb{R}^{n}$. A property of level sets of curves called the $\mathfrak{H}$-property will be defined, so-called because the level sets satisfying this property intersect like a subset of the hyperplanes in $\mathfrak{H}^{n}$. The level sets of $\{L(c) \mid c\in C\}$ will be said to have the $\mathfrak{H}$-property at $x\in \mathcal{T}_{g}$ if $T_{x}\mathcal{T}_{g}\simeq \mathbb{R}^{6g-6}$ has a basis \smash{$\bigl\{\frac{\partial}{\partial x_{1}}, \dots, \frac{\partial}{\partial x_{n}}\bigr\}$} that includes $\{\nabla L(c)(x) \mid c\in C\}$. Note that, when specifying the isomorphism between $T_{x}\mathcal{T}_{g}$ and $\mathbb{R}^{6g-6}$, a metric $\mathbf{g}$ on $\mathcal{T}_{g}$ is assumed, and the vector space isomorphism between $T_{x}\mathcal{T}_{g}$and $\mathbb{R}^{6g-6}$ is an isometry with respect to the standard metric of $\mathbb{R}^{6g-6}$. An equivalent definition of the $\mathfrak{H}$-property is that when the level sets of $\{L(c) \mid c\in C\}$ have the $\mathfrak{H}$-property at $x$, the following conditions are satisfied:
\begin{enumerate}\itemsep=0pt
\item[(1)] Each $C'\subset C$ determines exactly one face of $D(C,t,d)$ adjacent to $x$ contained in $L(C',t,\allowbreak d|_{C'})$. This will be called the face labelled by $C'$.
\item[(2)] The face from part 1 labelled by $C'$ has connected unit tangent cone.
\item[(3)] The tangent cone at $x$ of the face labelled by $C'\subset C$ has codimension equal to $|C'|-1$.
\end{enumerate}
Understanding the difference between linear independence of $\{\nabla L(c) \mid c\in C\}$ and level sets of~$\{L(c) \mid c\in C\}$ satisfying the $\mathfrak{H}$-property is slightly subtle, because the two properties are only different in dimensions at least 4. An example will be given later to illustrate this.

A tangent cone of a face of $D(C,t,d)$ labelled by $C'$ will be called \textit{full} if it has codimension~${|C'|-1}$, which is the codimension of $E(C', d|_{C'})$ when the lengths of curves in $C'$ have linearly independent gradients. A tangent cone to $D(C,t,d)$ is full if it has dimension equal to the dimension of $\mathcal{T}_{g}$.

\begin{Claim} At a $C$-generic point, the level sets of $\{L(c) \mid c\in C\}$ have the $\mathfrak{H}$-property.
\end{Claim}

The claim will be proven using twist paths. For each $c\in C$, at $x$, both $D(C\setminus \{c\},t,d)$ and~${L(C\setminus \{c\},t,d)\cap \partial D(C\setminus \{c\},t,d)}$ have a product structure. One term of the product (called an $\mathbb{R}$-factor below) is diffeomorphic to $\mathbb{R}$, and corresponds to the parameter on the $c^{*}$-twist paths. The intersection of the codimension 1 submanifold $L(c)^{-1}(t+d|_{c})$ with $L(C\setminus \{c\},t,d)\cap \partial D(C\setminus \{c\},t,d)$ cuts through this $\mathbb{R}$ factor; on one side, $c$ is longer than the other curves in $C$ and on the other it is shorter. The existence of an $\mathbb{R}$-factor corresponding to $c\in C$ guarantees the existence of a face labelled by $C\setminus \{c\}$; the face is the union of subintervals of~$\gamma_{c^{*},y}$ with~${y\in L(C,t,d)\cap \partial D(C,t,d)}$ on which $L(c)+d|_{c}<L(c_{i})+d|_{c_{i}}$ for all $c_{i}\in C\setminus \{c\}$. The $C$-generic assumption guarantees this is not empty, because there is a direction along $\gamma_{c^{*},y}$ for $y\in L(C,t,d)\cap \partial D(C,t,d)$ in which~$L(c)$ is decreasing and the lengths of the other curves in $C$ are stationary. For $c_{i}$ and $c_{j}\in C$, the face labelled by $C\setminus \{c_{i}, c_{j}\}$ lies between the faces labelled by $C\setminus \{c_{i}\}$ and $C\setminus \{c_{j}\}$. Faces labelled by smaller subsets are constructed similarly. Connectivity of unit tangent cones to faces is a consequence of linear independence, as is fullness of the tangent cones. This concludes the proof of the claim.

At all $C$-generic points that are not minima of $C$, the dimension of the tangent cone to the face labelled by $C'\subset C$ remains constant. This is a consequence of linear independence of~${\{\nabla L(c) \mid c\in C\}}$ at $C$-generic points.

The assumption that $x$ is $C$-generic will now be dropped. Suppose $x$ is a point of $L(C,t,d)$ on $\partial D(C,t,d)$ that is not a minimum. This assumption ensures that $D(C,t,d)$ has nonempty interior; $t$ can be decreased slightly and $L(C,t,d)\cap\partial D(C,t,d)$ does not become empty.

When $C$ does not fill, Proposition \ref{Thurstonprop} implies that at $x$, the tangent cone to $D(C,t,d)$ is full; these are directions in which the lengths of all the curves in $C$ are decreasing. When $C$ fills, $x$~is not in $\mathrm{Min}(C)$ as this would imply $x$ is a minimum of $C$. Hence by Definition \ref{defnmod} the tangent cone to $D(C,t,d)$ is also full. For any $x\in L(C,t,d) \cap \partial D(C, t, d)$ and any $c\in C$, Lemma \ref{LipschitzMap} gives a full tangent cone of the intersection of $ D(C\setminus \{c\}, t, d)$ and $L(c)_{\geq t+d|_{c}}$. This is the cone of directions in which $L(c)$ is increasing and the curves in $C\setminus \{c\}$ are decreasing. Every $c\in C$ labels a facet of~$D(C,t,d)$ with $x$ on the boundary. This is because, at $x$, the tangent cones of~$D(C,t,d)$ and of the intersection of $ D(C\setminus \{c\}, t, d)$ and $L(c)_{\geq t+d|_{c}}$ for every $c\in C$ are full, and by convexity the unit tangent cone of $ D(C\setminus \{c\}, t, d)$ is connected for every $c\in C$.

If the level sets $\{L(c) \mid c\in C\}$ have the $\mathfrak{H}$-property at $x\in \mathcal{T}_{g}$, this implies that the set of vectors $\{\nabla L(c)(x) \mid c\in C\}$ is linearly independent. Note that the converse is not true. One example is given by the gradients in $\mathbb{R}^{4}$ of the functions
\begin{equation*}
f_{1}=e^{x_{1}}+e^{-x_{2}},\qquad f_{2}=e^{x_{2}}+e^{-x_{3}},\qquad f_{3}=e^{-x_{3}}+e^{x_{4}},\qquad f_{4}=e^{-x_{4}}+e^{x_{1}}.
\end{equation*}
In this example, the gradients are linearly independent, except where $x_{2}=x_{4}$. Moreover, this set of functions has the property that at every point there is an open cone of directions in which all four functions are increasing, and for each $f_{i}$, $i=1,\dots, 4$, there is an open cone of directions in which $f_{i}$ is decreasing and the other functions increasing. Define $D(x)$ to be the intersection of the sublevel sets $\{f_{i}(x)\leq f_{i}(p) \mid i=1,2,3,4\}$. As for $D(C,t,d)$, there is a facet of $D(x)$ along a level set of each of the functions. This means that when $x_{2}=x_{4}$ at $x$, $\partial D(x)$ has no face labelled by $\{f_{1}, f_{2}\}$; the tangent space at $x$ to the submanifold on which $f_{1}=f_{2}$ consists of directions in which one of $\{f_{3}, f_{4}\}$ is increasing, and the other decreasing. Similarly for $f_{3}$ and~$f_{4}$. On a neighbourhood of the submanifold on which $x_{2}=x_{4}$, the level sets also cannot have the $\mathfrak{H}$-property.

As in the example, away from minima, Lemma \ref{LipschitzMap} and Proposition \ref{Thurstonprop} imply that if the level sets of $\{L(c) \mid c\in C\}$ do not have the $\mathfrak{H}$-property at $x$, this is also the case on a neighbourhood of~$x$. However, the set of $C$-nongeneric points has measure zero, and the level sets of $\{L(c) \mid {c\in C}\}$ have the $\mathfrak{H}$-property at $C$-generic points. This concludes the proof of the lemma.
\end{proof}

The key property of a minimal filling set of curves $C$ used to prove the previous lemma was that for every $c\in C$, $C\setminus\{c\}$ fills a strictly smaller subsurface of $\mathcal{S}_{g}$ than $C$. The next lemma is a generalisation of Lemma \ref{polyhedronargument}.

\begin{Lemma}
\label{neededtofill}
Let $C$ be a filling set of curves $($not necessarily minimal$)$ and suppose $c\in C$ has the property that $C\setminus\{c\}$ is not filling. Then $\nabla L(c)$ is linearly independent from $\{\nabla L(c') \mid c'\in C\setminus \{c\}\}$ everywhere in $\mathcal{T}_{g}$ away from minima of $C$.
\end{Lemma}

\begin{proof}
To prove this lemma, a ``relative version'' of the $\mathfrak{H}$-property will be defined. After it is shown that $C$ satisfies the $C\setminus\{c\}$-relative version of the $\mathfrak{H}$-property at points that are $\{c^{*}\}$-generic, the lemma follows from the same argument as in the proof of Lemma \ref{polyhedronargument}.

Suppose $x\in \mathcal{T}_{g}$ is not a minimum of $C$. Recall the three conditions given under the definition of the $\mathfrak{H}$-property. The curve $c$ will be said to satisfy the $C\setminus \{c\}$-relative $\mathfrak{H}$-property if
\begin{enumerate}\itemsep=0pt
\item[(1)] For every $C'\subset (C\setminus \{c\})$ that is the label of a face of $D(C\setminus\{c\},t,d(x))$ adjacent to $x$ contained in a level set of $E(C',t,d(x))$, there is exactly one face of $D(C,t,d(x))$ adjacent to $x$ contained in a level set of $E(C'\cup\{c\},t,d(x))$.
\item[(2)] The face from (1) labelled by $C'\cup\{c\}$ has unit tangent cone with the same number of connected components as the corresponding face of $D(C\setminus\{c\},t,d(x))$.
\item[(3)] For nonempty $C'$, the tangent cone at $x$ of the face of $D(C,t,d(x))$ labelled by $C'\cup \{c\}$ has dimension exactly one less that the corresponding face of $D(C\setminus\{c\},t,d(x))$ labelled by $C'$. When $C'$ is the empty set, $D(C,t,d(x))$ has a facet contained in a level set of $c$ and whose tangent cone has codimension 1 in $T_{x}\mathcal{T}_{g}$.
\end{enumerate}

At points that are $\{c^{*}\}$-generic, to prove that $c$ satisfies the $C\setminus \{c\}$-relative $\mathfrak{H}$-property, conditions (1), (2) and the first part of (3) are established using the same $\mathbb{R}$-factors as in the proof of the claim in Lemma \ref{polyhedronargument}. Lemma \ref{LipschitzMap} gives the existence of the facet in the second part of (3). As in the proof of Lemma \ref{polyhedronargument}, the assumption that $x$ is not a minimum of $C$ implies $x$ is not in $\overline{\mathrm{Min}(C)}$ hence the tangent cone to $D(C,t,d(x))$ is full. The tangent cone to~${D(C\setminus \{c\},t,d(x))\setminus D(C,t,d(x))}$ is full by Lemma \ref{LipschitzMap}. That the tangent cone to the facet labelled by $\{c\}$ has codimension 1 in $T_{x}\mathcal{T}_{g}$ is then a consequence of convexity of $D(C\setminus \{c\},t,d(x))$. This concludes the proof that $c$ satisfies the $C\setminus \{c\}$-relative $\mathfrak{H}$-property at $x$ when $x$ is $\{c^{*}\}$-generic.

As in the proof of Lemma \ref{polyhedronargument}, away from minima, $c$ satisfies the $C\setminus \{c\}$-relative $\mathfrak{H}$-property except on the complement of open sets. The lemma then follows from the same proof by contradiction as Lemma \ref{polyhedronargument}.
\end{proof}

The next lemma will be assumed from now on without comment.

\begin{Lemma}
\label{connectedminima}
If $C$ is a filling set of curves, whenever $E(C,d)$ is nonempty, $E(C,d)$ has a~minimum of $C$.
\end{Lemma}
\begin{proof}
Whenever $C$ fills, it contains a minimal filling set $C_{m}$. It was shown in Corollary~\ref{linear} that $E(C_{m}, d|_{C_{m}})$ has a minimum, and since $E(C,d)\subset E(C_{m}, d|_{C_{m}})$, the lengths of curves on~$E(C,d)$ cannot be less than the value at the minimum on $E(C_{m}, d|_{C_{m}})$. Since the intersection of $E(C,d)$ with the thick part of $\mathcal{T}_{g}$ is closed, the infimum is realised. As each connected component of $E(C,d)$ is closed, this also shows that each connected component of $E(C,d)$ has a minimum.
\end{proof}

\begin{Lemma}
\label{minimalem}
Let $C$ be a minimal filling set of curves and $C'\subsetneq C$. The only minima of $C$ occur in the intersection of $E(C,d)$ with $\mathrm{Min}(C)$, and $C'$ has no minima.
\end{Lemma}
\begin{proof}
It was shown in the proof of Lemma \ref{polyhedronargument} that for $c\in C$ there are no saddle points of~$L(c)|_{L(C\setminus\{c\},t,d)}$ and for $c\in C'$ there are no saddle points of $L(c)|_{L(C'\setminus\{c\},t,d)}$. It follows that each connected component of $E(C,d)$ and each connected component of $E(C',d)$ has at most one minimum of the restriction of the lengths of curves in $C$ respectively $C'$. By Lemma \ref{polyhedronargument}, the gradients of the lengths of curves in $C$, respectively $C'$, are linearly independent everywhere away from this minimum.

Every point $x$ in $\mathcal{T}_{g}$ is in $E(C,d)$ for some $d(x)$. It will now be shown that every $E(C,d)$ intersects $\mathrm{Min}(C)$ in a unique point.

That $E(C,d)$ intersects $\mathrm{Min}(C)$ in at most one point is an immediate consequence of the definition of $\mathrm{Min}(C)$. If $E(C,d)$ intersected $\mathrm{Min}(C)$ in two points, $p_{1}$ and $p_{2}$, then either the lengths of the curves in $C$ at both points are equal, or the lengths of the curves at one of the points, call it $p_{1}$, are smaller. In the first case, there are length functions with minima realised at both points, contradicting convexity. In the second case, a contradiction is obtained, because the length function with global minimum at $p_{2}$ is smaller at $p_{1}$.

Since $C$ is a minimal filling set of curves, by \cite[Lemma 9, Theorem 10]{SchmutzMorse}, $\mathrm{Min}(C)$ is a~continuously differentiable cell with empty boundary. As each $E(C,d)$ intersects $\mathrm{Min}(C)$ in at most one point, there is a differentiable bijection from $\mathrm{Min}(C)$ to a subset of elements of $\mathbb{R}^{|C|}$ whose entries sum to zero; a point $x\in \mathrm{Min}(C)$ is mapped to the tuple $d$ for which $E(C,d)$ passes through $x$.

\begin{Claim}
The set of loci $E(C,d)$ that intersect $\mathrm{Min}(C)$ with varying $d$ foliate $\mathcal{T}_{g}$.
\end{Claim}

To prove the claim, suppose $\{y_{i}\}$ is a sequence of points, each of which is in a stratum $E(C,d(y_{i}))$ intersecting $\mathrm{Min}(C)$, but with limit $y$ in a stratum $E(C,d(y))$ that does not intersect~$\mathrm{Min}(C)$. Let $\{x_{i}\}$ be a sequence of points in $\mathcal{T}_{g}$ converging to a point $x$, where $x_{i}$ is in the same locus as $y_{i}$ for every $i\in \mathbb{N}$. Since $d$ varies smoothly over $\mathcal{T}_{g}$, $\{d(y_{i})\}$ approaches $d(y)$ and~$\{d(x_{i})\}$ approaches $d(x)=d(y)$.

Note that there is a neighbourhood of $\mathrm{Min}(C)$ on which the claim holds. If the points of the sequence $\{x_{i}\}$ are chosen sufficiently close to $\mathrm{Min}(C)$, $x$ and every $x_{i}$ is in an $E(C,d(x_{i}))$ that intersects $\mathrm{Min}(C)$. Moreover, it is possible to choose the sequences in such a way that each $x_{i}$ and $y_{i}$ are arbitrarily close together, giving that $x$ and $y$ are in the same component of~$E(C,d(x))$. The claim then follows by contradiction. This concludes the proof of the lemma in the case of filling sets of curves.

For a nonfilling subset $C'$ of $C$, any minima of the lengths of the curves in $C'$ restricted to~$E(C',d)$ must occur within $\mathrm{Min}(C)$; otherwise this would give a minimum on some $E(C,d)$ away from $\mathrm{Min}(C)$. However, by \cite[Lemma 5]{SchmutzMorse}, the assumption that $C$ is a minimal filling set of curves implies that the gradients of the lengths of curves in $C'$ are everywhere linearly independent on $\mathrm{Min}(C)$. This rules out the existence of minima on $E(C',d)$, as these can only occur at a point at which $L(c)^{-1}(t+d|_{c})$ intersects $L(C'\setminus \{c\},t,d)$ nontransversely for any~${c\in C'}$.
\end{proof}

\begin{Corollary}
\label{linear}
Let $C$ be a minimal filling set of curves. When $E(C,d)$ is nonempty, it is a~connected, embedded submanifold of $\mathcal{T}_{g}$. The restriction to $E(C,d)$ of the lengths of the curves in $C$ is strictly convex and has a unique minimum.
\end{Corollary}

\begin{proof}
On a connected component of $E(C,d)$ that intersects $\mathrm{Min}(C)$, the gradients of the lengths of curves are linearly independent away from $\mathrm{Min}(C)$. Since $C=\{c_{1}, \dots, c_{k}\}$ is a minimal filling set of curves, it follows from \cite[Lemma 5]{SchmutzMorse} that $\bigl\{\nabla L(c_{2})-\nabla L(c_{1}), \dots, \nabla L(c_{k})-\nabla L(c_{1})\bigr\}$ is linearly independent everywhere on $\mathrm{Min}(C)$. Repeated applications of the pre-image lemma (see, for example, \cite[Corollary 5.14]{Lee}) then implies that each $E(C,d)$ passing through $\mathrm{Min}(C)$ is an embedded submanifold.

It was also was shown in the proof of Lemma \ref{minimalem} that on each connected component of~$E(C,d)$, the only nonregular point of the restriction of the lengths of curves in $C$ is a unique minimum where the connected component passes through $\mathrm{Min}(C)$ and that $E(C,d)$ can intersect $\mathrm{Min}(C)$ at most once. This proves convexity and connectivity.
\end{proof}

\subsection{Balanced loci}
\label{subbalanced}
This subsection discusses a property of loci. This concept is related to the way loci of filling sets of curves intersect corresponding sets of minima, leading to a trichotomy of three fundamentally different types of loci.

Suppose now that $C$ is any finite set of curves on $\mathcal{S}_{g}$ that intersect pairwise at most once. It follows from Lojasiewicz's theorem, \cite{Lojasiewicz1964} that $E(C,d)$ admits a triangulation, i.e., it is a simplicial complex. Consequently, $E(C,d)$ from equation~\eqref{withd} has a dimension, and in the interior of each simplex, it has a tangent space. In what follows, $v_{C}$ will always denote a vector or vector field that specifies a direction in which the lengths of the curves in $C$ are all increasing at the same rate. Whenever possible, $v_{C}$ is chosen to be contained in or on the boundary of the convex hull of $\{\nabla L(c) \mid c\in C\}$.

\begin{Definition}[$E(C,d)$ is balanced, semi-balanced or unbalanced at $x$]
Suppose $x$ is a point of $E(C,d)$. Then
\begin{itemize}\itemsep=0pt
\item{ $E(C,d)$ is \textit{balanced} at $x$ if $x$ is in $\mathrm{Min}(C)$ or there is a nonzero vector $v_{C}$ in $T_{x}E(C,d)$ in the interior of the convex hull of $\{\nabla L(c)(x) \mid c\in C\}$. Equivalently, if $V$ is the orthogonal complement of $v_{C}(x)$, $C$ is $V$-eutactic.}
\item{$E(C,d)$ is semi-balanced at $x$ if it is not balanced but $x$ is in $\mathrm{Min}(C')$ for $C'\subsetneq C$ or there is a nonzero vector $v_{C}$ on the boundary of the convex hull of $\{\nabla L(c)(x) \mid c\in C\}$.}
\item{$E(C,d)$ is unbalanced at $x$ if there exists a nonzero vector $v_{C}$ at $x$ and $E(C,d)$ is neither balanced nor semi-balanced.}
\end{itemize}
\end{Definition}

As will be explained later, the property of being balanced is independent of the choice of metric.

When $E(C,d)$ is balanced at $x$ and $x$ is not in $\mathrm{Min}(C)$, this implies that there is a direction~$v_{C}$ in $T_{x}E(C,d)$ in which the lengths of the curves in $C=\{c_{1}, \dots, c_{k}\}$ are increasing. While~${a_{1}\nabla L(c_{1})+\dots+a_{k}\nabla L(c_{k})}$, $a_{1}, \dots, a_{k}\in \mathbb{R}_{+}$ could be a direction in which some of the lengths of curves in $C$ are decreasing, when such a vector is tangent to $E(C,d)$, the constraints defining $E(C,d)$ force this to be a direction in which the lengths of all the curves in $C$ are increasing.

\begin{Remark}
When $C$ is a minimal filling set of curves, it was shown in Corollary \ref{linear} that~${E(C,d)\subset \mathcal{T}_{g}}$ is an embedded submanifold, in which case it has a tangent space at every point. The statement that $E(C,d)$ is balanced at every point therefore makes sense.
\end{Remark}
\begin{Proposition}
\label{balancedlem}
Let $C$ be a minimal filling set of curves. Then $E(C,d)$ is balanced at every point.
\end{Proposition}
\begin{proof}
The proof is essentially an application of Lagrange multipliers and the characterisation of $\mathrm{Min}(C)$ from \cite{SchmutzMorse}. Fix $C=\{c_{1}, \dots, c_{k}\}$. For $x$ not in $\mathrm{Min}(C)$, the locus $E(C,d)$ is balanced at $x$ if and only if is possible to find a $k$-tuple $(a_{1}, \dots, a_{k})$ in $\mathbb{R}_{+}^{k}$, with entries normalised to sum to one, and such that $v_{C}$ can be normalised to give $v_{C}=a_{1}\nabla L(c_{1})+\dots+a_{k}\nabla L(c_{k})$ in $T_{x}\mathcal{T}_{g}$.

\begin{Claim} For $x\in E(C,d)\setminus \mathrm{Min}(C)$, $v_{C}=a_{1}\nabla L(c_{1})+\dots+a_{k}\nabla L(c_{k})$ where $(a_{1}, \dots, a_{k})$ are the parameters of the point in $\mathrm{Min}(C)$ through which the locus $E(C,d)$ passes.
\end{Claim}

To prove the claim, recall that on every $E(C, d)$, by Corollary \ref{linear}, there is a unique point~$p(C,d)$ in the intersection $E(C,d)\cap\mathrm{Min}(C)$; the point $p(C,d)$ is the point at which the restrictions to $E(C,d)$ of the lengths of curves in $C$ have their minimum. Define
\begin{equation*}
 L(t):=\bigcup_{d}\{y\in E(C,d) \mid L(c)(y)=L(c)\bigl(p(C,d)\bigr)+t \, \forall c\in C\}.
\end{equation*}
Choose $t$ so that the connected component of $L(t)$ passes through the point $x$. Theorem 10 of \cite{SchmutzMorse} showed that $\mathrm{Min}(C)$ is a $C^{1}$-submanifold. Since $L(t)$ is the boundary of a tubular neighbourhood of a $C^{1}$-submanifold, it follows from a classical result, see \cite{Foote}, that $L(t)$ is also $C^{1}$. Lemma \ref{lemma14} implies that $\mathrm{Min}(C)$ has empty boundary.

Recall that the entries of the tuple $(a_{1}, \dots, a_{k})$ are normalised to sum to one. By construction, $x$ is a local minimum of the length function $a_{1}L(c_{1})+\dots+a_{k}L(c_{k})$ along any path in~$L(t)$ crossing~$E(C,d(x))$ transversely. It follows from the method of Lagrange multipliers that~${a_{1}\nabla L(c_{1})+\dots+a_{k}\nabla L(c_{k})}$ is normal to $L(t)$ everywhere in the intersection $E(C,d(x))\cap L(t)$, and hence parallel to the normal vector field $\frac{\partial}{\partial t}$ on $L(t)$. In $E(C,d(x))$, the common gradient of the length functions, a vector field parallel to $v_{C}$, is also always orthogonal to the level sets of the lengths of curves in $C$. This proves the claim.

As this construction works for any value of $t>0$, it follows that the $k$-tuple $d(x)$ defining the locus $E(C,d(x))$ passing through $x$ is the value of $d|_{C}$ for which $E(C, d(x))$ passes through the point in $\mathrm{Min}(C)$ representing the minimum of the length function $a_{1}L(c_{1})+\dots+a_{k} L(c_{k})$, and $v_{C}$ is parallel to $a_{1}\nabla L(c_{1})+\dots+a_{k}\nabla L(c_{k})$ everywhere on $E\big(C,d(x)\big)$. This concludes the proof.
\end{proof}

\begin{Proposition}
\label{balancedlemnonfilling}
Let $C'$ be a subset of a minimal filling set, or a multicurve possibly containing separating curves. Then $E(C')$ is balanced at every point.
\end{Proposition}

\begin{Remark}
\label{balancednearby}
Unlike the case for minimal filling sets, $E(C',d)$ is not balanced for every possible $d$. However, since $\{\nabla L(c)(x) \mid c\in C'\}$ is linearly independent at $x\in E(C')$, this linear independence holds on a neighbourhood of $x$ in $\mathcal{T}_{g}$, as does the existence of a $v_{C'}$ in the convex hull of $\{\nabla L(c) \mid c\in C'\}$. It follows that points in loci $E(C',d)$ sufficiently close to $E(C')$ are also balanced.
\end{Remark}

\begin{proof}
When $C'$ is a subset of a minimal filling set, it follows from Lemmas \ref{polyhedronargument} and \ref{minimalem} that the loci $E(C',d)$ are embedded submanifolds that foliate $\mathcal{T}_{g}$ with varying $d$. When $C'$ is a~multicurve, this follows from the fact that the lengths of curves in $C'$ determine a subset of Fenchel--Nielsen coordinates defined everywhere on $\mathcal{T}_{g}$. The proof of Proposition \ref{balancedlem} therefore generalises as soon as a replacement for the set of minima $\mathrm{Min}(C)$ can be found.

By Proposition \ref{Thurstonprop}, $\mathrm{Min}(C')$ is empty in $\mathcal{T}_{g}$. However, nonfilling subsets can be understood to have sets of minima in the metric completion $\overline{\mathcal{T}}_{g}$ of $\mathcal{T}_{g}$ with respect to the Weil--Petersson metric. Points on the boundary of $\mathcal{T}_{g}$ with respect to this completion can be understood to be noded surfaces that have been pinched along multicurves, see, for example, \cite{wolpert2005}. In \cite{SchmutzMorse}, sets of minima were also defined for noded surfaces. It follows from Riera's formula~\cite{Riera} that a~necessary condition for a length function $L(A, C')$ to have a minimum is that the multicurve consisting of the boundary of the subsurface filled by $C'$ has been pinched. Within such noded surfaces, $C'$ is a minimal filling set, and it is possible to define a set of minima that will be denoted by~$\mathrm{Min}(C')^{\infty}$. Not every $E(C',d)$ will necessarily have a point in $\mathrm{Min}(C')^{\infty}$ in its completion; the proof consequently reduces to the problem of showing that $E(C')$ does.

\begin{Corollary}[{\cite[Corollary 21]{wolpert2005}}]
\label{wolcor}
Denote by $\mathrm{Sys}(m)^{\infty}$ a stratum in $\overline{\mathcal{T}}_{g}\setminus \mathcal{T}_{g}$ corresponding to noded surfaces on which the multicurve $m$ has been pinched. Denote by $L(m)$ the length function~$\sum_{c\in m}L(c)$.
The distance $d(x, \mathrm{Sys}(m)^{\infty})$ of a point $x\in \mathcal{T}_{g}$ to $\mathrm{Sys}(m)^{\infty}$ has the expansion
\begin{equation*}
d(x, \mathrm{Sys}(m)^{\infty})=\sqrt{2\pi L(m)}+\mathcal{O}\bigl(L(m)^{2}\bigr).
\end{equation*}
\end{Corollary}

When $C'$ is a multicurve $m$, $\mathrm{Min}(C')^{\infty}$ is contained in $\mathrm{Sys}(m)^{\infty}$, which is in the metric completion of $E(C')$. This proves the proposition for multicurves.

When $C'$ is not a multicurve, denote by $m$ a multicurve contained in $C'$ with largest cardinality, and $c_{1}$ a curve in $C'\setminus m$. Then by the collar lemma, $L(c_{1})$ is zero on $\mathrm{Sys}(\{c_{1}\})^{\infty}$ and infinite on $\mathrm{Sys}(m)^{\infty}$, similarly $L(m)$ is zero on $\mathrm{Sys}(m)^{\infty}$ and infinite on $\mathrm{Sys}(\{c_{1}\})^{\infty}$. Again by Riera's formula, $\mathrm{Min}(m\cup\{c_{1}\})^{\infty}$ must be in a noded surface pinched along the multicurve on the boundary of the subsurface filled by $m\cup\{c_{1}\}$. This multicurve is disjoint from both $m$ and~$c_{1}$, so this analogue of a set of minima has the analogues of sets of minima of subsets of~${m\cup\{c_{1}\}}$ on its boundary, generalising Lemma \ref{lemma14}.

Let $\gamma$ be a path in the thin part of $\mathcal{T}_{g}$ near $\mathrm{Min}(m\cup\{c_{1}\})^{\infty}$ with one endpoint close to~$\mathrm{Sys}(\{c_{1}\})^{\infty}$ and the other close to $\mathrm{Sys}(m)^{\infty}$. Using the estimate in Corollary \ref{wolcor}, it follows from the intermediate value theorem applied to $\gamma$ that $E(m\cup\{c_{1}\})$ has a point in $\mathrm{Min}(m\cup\{c_{1}\})^{\infty}$ in its metric completion. This proves the proposition when $C'=m\cup \{c_{1}\}$. Now let $c_{2}$ be a~curve in $C'\setminus (m\cup \{c_{1}\})$ chosen such that $c_{1}$ intersects a curve in $m\cup \{c_{1}\}$; this is always possible because $m$ was assumed to be a multicurve of maximal cardinality. The same argument proves the proposition for $C'=m\cup \{c_{1},c_{2}\}$. This construction is then iterated further if necessary.
\end{proof}

\begin{Remark}
\label{subsetrem}
When using Proposition \ref{balancedlemnonfilling}, it can be helpful to know that any pair $\{c_{1}, c_{2}\}$ of curves intersecting exactly once is contained in a minimal filling set. This follows from the observations that the action of $\Gamma_{g}$ on curves takes minimal filling sets to minimal filling sets, $\Gamma_{g}$ acts transitively on pairs of nonseparating curves intersecting exactly once, and any minimal filling set contains a pair of curves that intersect once.
\end{Remark}

\begin{Corollary}
\label{polytopecor}
A stratum $\mathrm{Sys}(C)$ is on the boundary of $\mathrm{Sys}(\{c\})$ for every $c\in C$.
\end{Corollary}
\begin{proof}
Let $x$ be a point in $\mathrm{Sys}(C)$. It follows from Proposition \ref{balancedlemnonfilling} and Remark \ref{subsetrem} that for every $c_{1}, c_{2}\in C$, $\mathbf{g}(\nabla L(c_{2})(x),\nabla L(c_{1})(x))<\mathbf{g}(\nabla L(c_{1})(x), \nabla L(c_{1})(x))$. For every $c\in C$, it follows from local finiteness that the vector $-\nabla L(c)(x)$ is in the tangent cone of $\mathrm{Sys}(\{c\})$.
\end{proof}

\begin{Proposition}
The property of being balanced at $x$ is independent of the choice of metric.
\end{Proposition}
\begin{proof}
The proof is basically a repeat of the arguments given in Propositions \ref{balancedlem} and \ref{balancedlemnonfilling}. Suppose $E(C,d)$ is balanced at $x$ with respect to a metric $\mathbf{g}$. Setting $C=\{c_{1}, \dots, c_{k}\}$, one obtains that $v_{C}$ is proportional to $\sum_{i=1}^{k} a_{i}\nabla L(c_{i})$ where $a_{i}\in \mathbb{R}_{+}$, $i=1, \dots, k$ and $\sum_{i=1}^{k}a_{i}=1$. Define $L:=\sum_{i=1}^{k} a_{i}L(c_{i})$. The proposition will be proven by showing that $\nabla L(x)$ is proportional to $v_{C}(x)$ for any metric.

The functions $\mathcal{T}_{g}\rightarrow \mathbb{R}_{+}$ obtained by taking components of the $k$-tuple $d$ will be denoted by $d_{1}, \dots, d_{k}$. On a neighbourhood of $x$, it is possible to obtain a subset of independent coordinates $\{d'_{1}, \dots, d'_{n}\}$ with the property that each $d'_{i}$ is a linear combination of $d_{1}, \dots, d_{k}$, and~${\mathrm{Span}\{\nabla d'_{1}, \dots, \nabla d'_{n}\}=\mathrm{Span}\{\nabla d_{1}, \dots, \nabla d_{k}\}}$.

The vectors $v_{C}(x)$ and a spanning set for $\{\nabla L(c)(x) \mid c\in C\}^{\perp}$ in $T_{x}\mathcal{T}_{g}$ span the subspace $\{\nabla d_{1}, \dots, \nabla d_{k}\}^{\perp}$ of $T_{x}\mathcal{T}_{g}$ on which the tuple $d$ is stationary. Note that this property is independent of $\mathbf{g}$, as the subspace of $T_{x}\mathcal{T}_{g}$ on which $d$ is stationary can be defined in terms of tangents to intersections of level sets, as can $\{\nabla L(c)(x) \mid c\in C\}^{\perp}$. Taking $d'_{1}, \dots, d'_{n}, L$ to be a~subset of independent coordinates on a neighbourhood of $x$ in $\mathcal{T}_{g}$, calculating differentials in this coordinate chart gives that $dL(x)$ is contained in \smash{$\bigl(\{\nabla d'_{1}, \dots, \nabla d'_{n}\}^{\perp}\bigr)^{\flat}$}. Also by construction, \smash{$\nabla L(x)=\sum_{i=1}^{k} a_{i}\nabla L(c_{i})(x)$}. Consequently, $dL(x)=v_{C}^{\flat}(x)$, as required.
\end{proof}

\begin{figure}[ht]
\centering
\includegraphics[width=0.76\textwidth]{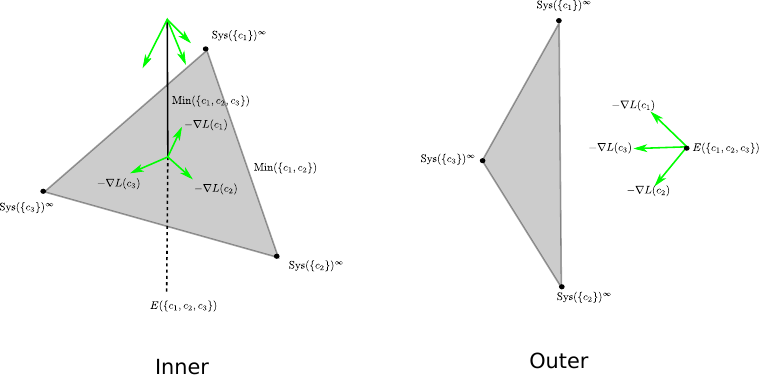}
\caption{Inner versus outer loci. In the figure on the right, the locus should be regarded as coming out of the page. A filling set of curves that intersect pairwise at most once must have cardinality at least four. This figure represents a lower-dimensional analogue, where $\mathrm{Sys}(\{c_{i}\})^{\infty}$ denotes a point in the completion of $\mathcal{T}_{g}$ with respect to the Weil--Petersson metric at which the curve $c_{i}$ has been pinched.}
\label{innerouter}
\end{figure}

Once again allowing $C$ to be any filling set of curves that intersects pairwise at most once, now that $\mathrm{Min}(C)$ can have a boundary, a new phenomenon occurs.

\textbf{Inner, outer and borderline loci.} A connected component of $E(C,d)$ is
\begin{itemize}\itemsep=0pt
\item\textit{inner} if it passes through a point of $\mathrm{Min}(C)$,
\item\textit{borderline} if it passes through a point of $\partial \mathrm{Min}(C)$,
\item\textit{outer} if it is disjoint from $\overline{\mathrm{Min}(C)}$.
\end{itemize}
Lower-dimensional analogues of inner and outer loci are shown in Figure~\ref{innerouter}. As $E(C,d)$ can intersect $\overline{\mathrm{Min}(C)}$ in at most one point, each locus is exactly one of these three types.

An inner locus $E(C, d)$ has the property that at the global minimum point of the lengths of curves in $C$ restricted to $E(C, d)$, there is no direction in which the lengths of all the curves in~$C$ are decreasing, whereas there is an open cone of such directions at the minimum point when~$E(C, d)$ is outer.

The loci $E(C,d)$ studied here so far were all embedded submanifolds of dimension greater than one that foliate $\mathcal{T}_{g}$ for varying $d$. The arguments of Proposition \ref{balancedlem} suggest that for filling~$C$, such loci are
\begin{itemize}\itemsep=0pt
\item balanced if they are inner
\item semi-balanced if they are borderline
\item unbalanced if they are outer.
\end{itemize}
It would be interesting to know to what extent this generalises.

For a filling set $C$, if $\mathrm{Sys}(C)\subset E(C)$ is not balanced at some point $x$, then by Proposition~\ref{balancedlem}, $C$ strictly contains a minimal filling set of curves. The next proposition can be used to characterise a filling subset of $C$.

\begin{Proposition}
\label{plus1}
Suppose $C$ is a filling set of curves, and at a point $x$ that is not a minimum of $C$, a nonzero $v_{C}(x)$ exists which cannot be chosen to be in the interior of the convex hull of $\{\nabla L(c)(x) \mid c\in C\}$. Then the vertices of the face $f$ of the convex hull of $\{\nabla L(c)(x) \mid c\in C\}$ closest to or containing $v_{C}(x)$ are the gradients of a filling set of curves $c(f)$.
\end{Proposition}

\begin{Remark}
\label{balancedrem}
The set of curves $c(f)$ from Proposition \ref{plus1} has the following property in common with a minimal filling set of curves: A nonzero vector $v_{c(f)}(x)$ exists in the interior of the convex hull of $\{\nabla L(c)(x) \mid c\in c(f)\}$. Identifying $T_{x}\mathcal{T}_{g}$ with $\mathbb{R}^{6g-6}$ by a vector space isomorphism that preserves the inner product given by $\mathbf{g}$, the vector $v_{c(f)}(x)$ is obtained by adding to $v_{C}(x)$ a multiple of the normal vector to $f$ representing the shortest path between the vertex of $v_{C}(x)$ and the convex hull.
\end{Remark}

\begin{proof}[Proof of Proposition~\ref{plus1}]
The motivation for the proof is to compare the assumption on the gradients with the statement of Lemma \ref{lemma14}. In the projection to the orthogonal complement of $v_{c(f)}(x)$ in $T_{x}\mathcal{T}_{g}$, the gradients $\{\nabla L(c)(x) \mid c\in C\}$ look like the gradients at a point in \smash{$\mathrm{Min}(c(f))\subset\overline{\mathrm{Min}(C)}$}, where~$\mathrm{Min}(c(f))$ is only nonempty if $c(f)$ fills. The intuition is that $c(f)$ corresponds to the set of minima on the boundary of $\mathrm{Min}(C)$ closest to $x$; such a boundary point is reached by decreasing the lengths of curves in $C$.

For every $c\in (C\setminus c(f))$, the inequality
\begin{equation}
\label{ineq}
\mathbf{g}(\nabla L(c), v_{c(f)})\geq\mathbf{g}(v_{c(f)}, v_{c(f)})
\end{equation}
is assumed in the statement of the proposition to hold at $x$, where the assumption is made that~$v_{c(f)}$ is normalised in such a way that it represents a point in $\mathbb{R}^{6g-6}$ contained in the polytope with vertices given by $\{\nabla L(c)\mid c\in c(f)\}$.

Denote by $C_{\mathrm{ineq}}$ the connected component containing $x$ of the subset of $E(c(f),d(x))$ on which inequality~\eqref{ineq} holds. A point at which $v_{c(f)}=0$ is assumed to be in $C_{\mathrm{ineq}}$ if it is a limit point of a~sequence on which inequality~\eqref{ineq} holds.

The argument is greatly simplified at a number of points by arguing that pathological cases (usually to do with linear dependence of gradients) can only occur when one or more curves can be discarded from the set without changing the subsurface of $\mathcal{S}_{g}$ filled by the set. To begin with, it will be shown that the following assumptions can be made without loss of generality.
\begin{enumerate}\itemsep=0pt
\item[(1)] The set of vectors $\{\nabla L(c) \mid c\in (C\setminus c(f))\}$ is linearly independent in $\mathcal{T}_{g}$ away from minima of $C\setminus c(f)$, and each one is linearly independent from $\{\nabla L(c) \mid c\in c(f)\}$ away from minima of $C$.
\item[(2)] The set of vectors $\{\nabla L(c)(x) \mid c\in c(f)\}$ is linearly independent everywhere in $C_{\mathrm{ineq}}$ where $v_{c(f)}$ is nonzero. A nonzero vector $v_{c(f)}$ whose existence is guaranteed by this linear independence can be chosen to be contained in the convex hull of $\{\nabla L(c) \mid c\in c(f)\}$ and is in the tangent cone to $E(c(f),d(x))$.
\item[(3)] Strict inequality in inequality~\eqref{ineq} holds at $x$.
\end{enumerate}
To see that assumption (1) can be made without loss of generality, suppose some of the curves~${C\setminus c(f)}$ have lengths with linearly dependent gradients at a point that is not a minimum of~$C\setminus c(f)$. It follows from Lemma \ref{neededtofill} that the proof of the proposition can be simplified by removing these curves from $C$. Doing so does not destroy the property that inequality~\eqref{ineq} holds for the remaining curves in $C\setminus c(f)$. Similarly, if there is a curve in $C\setminus c(f)$ for which the gradient of the length is linearly dependent on $\{\nabla L(c) \mid c\in c(f)\}$ away from minima of $C$, discarding this curve from $C$ still results in a filling set. If all the curves $C\setminus c(f)$ have lengths with gradients linearly dependent on $\{\nabla L(c) \mid c\in c(f)\}$ at a point that is not a minimum of $C$, it follows that $c(f)$ fills as required. This concludes the proof that the first assumption can be made without loss of generality.

\begin{figure}[ht]
\centering
\includegraphics[width=0.7\textwidth]{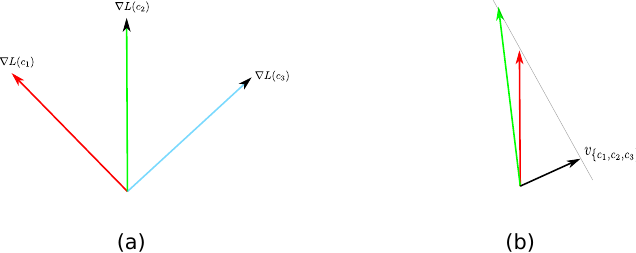}
\caption{If no nonzero $v_{c(f)}(y)$ exists because the dimension of $\mathrm{Span}(\{\nabla L(c) \mid c\in c(f)\})$ drops at $y$ (part~(a)), then at points near $y$ where a nonzero $v_{c(f)}$ exists, this is not in the convex hull of $\{\nabla L(c) \mid c\in c(f)\}$ (part (b)), viewed from the side. }
\label{tipover}
\end{figure}

At every point $y$ in $C_{\mathrm{ineq}}$ at which $v_{c(f)}$ is nonzero, it is possible to assume without loss of generality that the set of vectors $\{\nabla L(c')(y) \mid c'\in c(f)\}$ is linearly independent. Otherwise, by Lemma \ref{neededtofill} there is a subset $c'(f)$ of $c(f)$ that fills the same subset of $\mathcal{S}_{g}$ as $c(f)$ and for which~${v_{c'(f)}(y)=v_{c(f)}(y)}$. It is then possible to prove the proposition with $c'(f)$ is place of~$c(f)$. The same argument shows that it is possible to assume without loss of generality that if a nonzero~$v_{c(f)}(y)$ exists for $y\in C_{\mathrm{ineq}}$, it is in the tangent cone to $E(c(f),d(x))$ at $y$. It is also possible to assume without loss of generality that $v_{c(f)}$ is contained in the convex hull of~$\{\nabla L(c') \mid c'\in c(f)\}$ at every point in $C_{\mathrm{ineq}}$ where $v_{c(f)}$ is nonzero because otherwise there is a point $y'$ satisfying the conditions of the proposition with $c'(f)\subset c(f)$ in place of $c(f)$. This concludes the proof that the second assumption can be made without loss of generality at points of $C_{\mathrm{ineq}}$ where a nonzero $v_{c(f)}$ exists.

Strict convexity of length functions with respect to the Weil--Petersson metric, proven in \cite{Wolpert}, implies that if equality holds in inequality~\eqref{ineq} for a curve $c\in (C\setminus c(f))$, shifting $x$ an arbitrarily small amount in the direction of $\nabla L(c)(x)$ gives a point $x'$ at which strict inequality holds. By assumption (1), it is then possible to work in the locus containing $x'$ instead of $E(C,d(x))$, showing that assumption (3) can also be made without loss of generality.

To understand what happens when $v_{c(f)}(y)=0$ for $y\in C_{\mathrm{ineq}}$, suppose $y$ is on the boundary of the subset of $C_{\mathrm{ineq}}$ on which $v_{c(f)}$ is zero. It follows from assumption (2) applied to points of $C_{\mathrm{ineq}}$ at which $v_{c(f)}$ is nonzero that~$y$ is either in \smash{$\overline{\mathrm{Min}(c(f))}$} or is a point at which linear independence of $\{\nabla L(c') \mid c'\in c(f)\}$ breaks down. The former implies $c(f)$ fills, so suppose from now on that it is the latter. Definition~\ref{defnmod} implies that there is an open cone of directions at~$y$ in which the lengths of all the curves in $c(f)$ are decreasing. By assumption, $y$ is a limit point of a sequence on which linear independence of $\{\nabla L(c') \mid c'\in c(f)\}$ holds and on which there is a nonzero $v_{c(f)}$ in the convex hull of $\{\nabla L(c') \mid c'\in c(f)\}$. As illustrated in Figure~\ref{tipover}, this prevents linear independence of $\{\nabla L(c') \mid c'\in c(f)\}$ from breaking down at $y$, because for all points sufficiently close to $y$ at which linear independence holds, $v_{c(f)}$ could not be in the convex hull of $\{\nabla L(c') \mid c'\in c(f)\}$. This shows that $v_{c(f)}=0$ on $C_{\mathrm{ineq}}$ implies that $c(f)$ fills. It will therefore be assumed that $v_{c(f)}$ is nonzero on $C_{\mathrm{ineq}}$.

Next, assume to begin with that $|C\setminus c(f)|=1$. The locus $E(C,d(x))$ is a separating subset of $E(c(f),d(x))$. One side of $E(c(f),d(x))\setminus E(C,d(x))$ contains the loci in which the curve $c$ in~$C\setminus c(f)$ is longer than the curves in $c(f)$, and the other contains the loci on which it is shorter. The loci $E(C,d)$ for varying $d$ contained in $E(c(f), d(x))$ are shown schematically in Figure~\ref{q1}, as well as other objects that will be defined in this proof.

\begin{figure}[ht]
\centering
\includegraphics[width=0.6\textwidth]{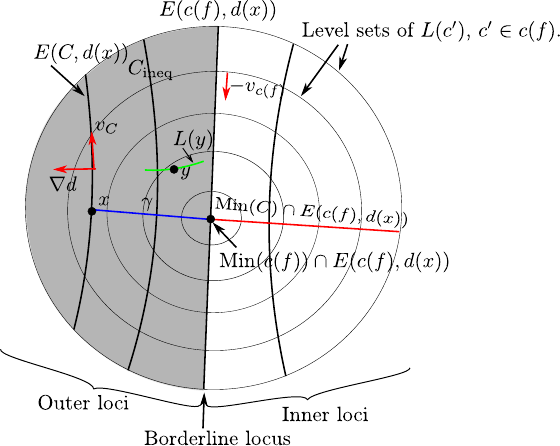}
\caption{A schematic drawing of $E(c(f),d(x)$. The subset $C_{\mathrm{ineq}}$ of $E(c(f),d(x))$ is shaded. The concentric circles are level sets of $L(c')$ for $c'\in c(f)$, and the curved vertical lines are loci $E(C,d)$ for different values of $d$. }
\label{q1}
\end{figure}

In what follows, within $E(c(f),d(x))$ it will be convenient to confuse the tuple $d$ with a~function $d\colon E(c(f),d(x))\rightarrow \mathbb{R}$, $x\mapsto L(c)(x)-L(c')(x)$ for a choice (it does not matter which) of curve $c'\in c(f)$.

\begin{Claim}
 Under the assumptions listed above, if strict inequality in inequality~\eqref{ineq} holds at $x$, it holds everywhere on the connected component $E'(C,d(x))$ of $E(C,d(x))$ containing $x$.
 \end{Claim}

The claim will be proven by contradiction. If the claim were false, it follows from the linear independence assumptions that there is a point $y_0$ on $E(C,d(x))$ at which equality in inequality~\eqref{ineq} holds. Choose a point $y$ on $E'(C,d(x))$ close to $y_0$ at which inequality in inequality~\eqref{ineq} holds, and for which $y$ is not a minimum of $C$. Setting $C=\{c_{1}, \dots, c_{k}\}$, $c=c_{1}$, there is a~$b_{1}<0$, $b_{i}>0$ for $i=2, \dots,k$ with $\sum_{i=1}^{k}b_{i}=1$ such that $v_{C}(y)=\sum_{i=1}^{k}b_{i}\nabla L(c_{i})(y)\neq 0$. The point~$y$ is chosen such that $|b_{1}|$ is sufficiently small that the function $\sum_{i=1}^{k}b_{i}L(c_{i})$ is convex on a neighbourhood of~$y$ containing $y_{0}$.

An argument similar to the proof of Proposition \ref{balancedlem} will now be given, but using an analogue of $\mathrm{Conv}(C)$ (defined in Section \ref{SchmutzSection}) in place of $\mathrm{Min}(C)$. The general idea is to use the linear independence assumptions and convexity to derive a contradiction to the existence of $y_{0}$.

Within some neighbourhood of $y$, the flow-box theorem guarantees the existence of a continuously differentiable codimension 1 submanifold $L(y)$ embedded in $E(c(f),d(x))$ orthogonal to the vector field $v_{C}=v_{\{c\}\cup c(f)}$ and passing through $y$. On a neighbourhood of $y$, a locus $E(C,d)$ near $E'(C,d(x))$ intersects $L(y)$ along a level set of $L(C)|_{E(C,d)}$. Define
\begin{align*}
L(y,t)=\bigcup_{d\in \mathcal{N}(d(x))}\{{}&z\in E(C,d) \mid L(c_{i})(z)=t+L(c_{i})(z_{0}) \, \forall c_{i}\in C, \\
&\text{where }z_{0}\in E(C,d)\cap L(y)\},
\end{align*}
where $t$ takes values in a neighbourhood of $0$ in $\mathbb{R}$, and $\mathcal{N}(d(y))$ is a neighbourhood of $d(x)=d(y)$ in $\mathbb{R}$.

It follows from convexity of the length functions $L(c)$ and $L(c_{2})$ and local convexity of $\sum_{i=1}^{k}b_{i}L(c_{i})$ that $E'(C,d(x))\cap L(y)$ is the minimum of the restriction of the function $\sum_{i=1}^{k}b_{i}L(c_{i})$ to $L(y)$. Similarly, $E'(C,d(x))\cap L(y,t)$ is the minimum of the restriction of the function $\sum_{i=1}^{k}b_{i}L(c_{i})$ to $L(y,t)$. The method of Lagrange multipliers then shows that on a neighbourhood of $y$ containing $y_{0}$, \[
v_{C}|_{E'(C,d(x))}=\sum_{i=1}^{k}b_{i}\nabla L(c_{i})\neq 0,
\]
 contradicting the assumption that equality in inequality~\eqref{ineq} occurs at $y_{0}$. This concludes the proof of the claim.

The minima of $C$ on $E'(C,d(x))$ are a connected set, because the linear independence assumptions rule out any nonregular points of $L(C)|_{E'(C,d(x))}$ that are not minima. Within $C_{\mathrm{ineq}}$ the loci $E(C,d)$ for varying $d$ make up the leaves of a codimension 1 foliation. This is because every connected component of $E(C,d)$ in $C_{\mathrm{ineq}}$ is noncompact, as otherwise there would be loci in $C_{\mathrm{ineq}}$ with maxima of $C$, contradicting the linear independence assumptions. Since $C$ fills, level sets of $L(C)|_{E(C,d)}$ are compact. Consequently, there is no connected component of a locus~$E(C,d)$ in $C_{\mathrm{ineq}}$ contained in a level set of $L(C)$. The function $d$ can only take regular values in the interior of $C_{\mathrm{ineq}}$ because otherwise there would be a nonregular value $d^{*}$ of $d$ such that assumption (1) breaks down at every point of a connected component $E'(C,d^{*})$ of $E(C,d^{*})$ in~$C_{\mathrm{ineq}}$ that is not a minimum of $C$.

\begin{Remark}
\label{holdsonnbhd}
Since assumptions (1) and (2) hold on a neighbourhood of $C_{\mathrm{ineq}}$, and $v_{c_{f}}$ is assumed to be nonzero in $C_{\mathrm{ineq}}$, it also follows that there is a neighbourhood of $C_{\mathrm{ineq}}$ in~$E(c(f),d(x))$ on which the function $d$ can only take regular values.
\end{Remark}

Let $\gamma$ be a continuous path traced out by minima of $C$ on $E(c(f),d(x))$, starting with the minimum of $C$ in $E(C,d(x)) \subset E(c(f),d(x))$ and passing through the loci in the direction of~$-v_{c(f)}$. Since $d$ only takes regular values on $C_{\mathrm{ineq}}$, strict inequality in inequality~\eqref{ineq} holds along the intersection of $\gamma$ with the interior of $C_{\mathrm{ineq}}$. Moreover, by assumption (2) and Remark~\ref{holdsonnbhd}, a~point at which $\gamma$ ends or leaves $C_{\mathrm{ineq}}$ is a point at which $v_{c(f)}=0$. The proposition will be proven by showing that $\gamma$ contained in $C_{\mathrm{ineq}}$ leads to a contradiction.

If $\gamma$ is contained in $C_{\mathrm{ineq}}$, $\gamma$ approaches a noded surface in $\overline{\mathcal{T}}_{g}$ on which a multicurve has been pinched. For example, when $d(x)=0$, as discussed in the proof of Lemma \ref{balancedlemnonfilling}, this is the noded surface on which the multicurve on the boundary of the subsurface filled by $c(f)$ has been pinched. Since $C$ fills, this multicurve intersects a curve in $C$. It follows from the collar lemma that $\gamma$ cannot be a path along which the lengths of curves in $C$ decrease monotonically, giving the required contradiction.

When $|C|>|c(f)|+1$, assumption (1) makes it possible to replace $L(c)$ in the previous argument by the length function
\begin{equation*}
L(C\setminus c(f)):=\frac{1}{|C\setminus c(f)|}\sum_{c'\in (C\setminus c(f))}L(c').
\end{equation*}
Along $\gamma$, the fact that the length function $L(C\setminus c(f))$ is decreasing faster in the direction of~$v_{c(f)}$ than the lengths of curves in $c(f)$ does not mean that the length of every curve in~${C\setminus c(f)}$ is decreasing faster in the direction of $-v_{c(f)}$ than the lengths of curves in $c(f)$. If, for example, $C\setminus c(f)=\{c_{1}, c_{2}\}$, and at the point $p$ at which $\gamma$ leaves $C_{\mathrm{ineq}}$ the inequality~\eqref{ineq} holds for~$c_{1}$ but not~$c_{2}$, then $\gamma$ passes through a locus on which the assumptions of the proposition hold with $c(f)$ replaced by $c(f)\cup\{c_{2}\}\subset C$. The argument then shows $c(f)\cup \{c_{2}\}\subset C$ fills, and can subsequently be applied with $c(f)\cup \{c_{2}\}$ in place of $C$ to show that $c(f)$ fills. This argument can be generalised to the case with $|C\setminus c(f)|>2$.
\end{proof}

\section{Proof of the Morse--Smale properties}
\label{mainproof}

The purpose of this section and the next is to show that $\mathcal{P}_{g}$ contains a choice of unstable manifolds of the critical points of $f_{\mathrm{sys}}$. This section studies the behaviour of $\mathcal{P}_{g}$ away from critical points and boundary points of $f_{\mathrm{sys}}$, i.e., at points at which the cone of increase is full. It will be shown that the locally top-dimensional strata making up the cells of $\mathcal{P}_{g}$ are balanced, and that restricting $f_{\mathrm{sys}}$ to $\mathcal{P}_{g}$ does not create critical points. In this way, the intuition of Thurston that $f_{\mathrm{sys}}$ is increasing towards $\mathcal{P}_{g}$ is encoded.

The main results of this section will be proven by generalising the properties of strata of minimal filling sets to the locally top-dimensional strata of $\mathcal{P}_{g}$. The simple case is described in the next lemma.

\begin{Lemma}
\label{linearlyindependent}
Suppose $x\in \mathrm{Sys}(C)$ and the set of vectors $\{\nabla L(c)(x) \mid c\in C\}$ is linearly independent. Then $\mathrm{Sys}(C)$ is on the boundary of $\mathrm{Sys}(C')$ for every nonempty $C'\subsetneq C$.
\end{Lemma}
\begin{proof}
For any fixed $C'\subset C$, $E(C)\subset E(C')$, so $E(C')$ is not empty. The linear independence assumption then ensures that the codimension of $E(C')$ in $\mathcal{T}_{g}$ is equal to $|C'|-1$, and that the tangent space at $x$ to the level set of $E(C')$ passing through $x$ is given by the orthogonal complement of $\mathrm{Span}\{\nabla L(c) \mid c\in C'\}$. Choose $c_{1}\in C$. It follows from local finiteness that the nonzero vector obtained by projecting $\nabla L(c_{1})(x)$ onto the orthogonal complement of $\mathrm{Span}\{\nabla L(c) \mid c\in C\setminus \{c_{1}\}\}$ is in the tangent cone to $\mathrm{Sys}(C\setminus\{c_{1}\})$. At a point $y\in \mathrm{Sys}(C\setminus\{c_{1}\})$ that can be made arbitrarily close to $x$, for any choice of $c_{2}\in (C\setminus\{c_{1}\})$, a stratum $\mathrm{Sys}(C\setminus\{c_{1},c_{2}\})$ adjacent to $y$, and hence also to $x$, can be constructed similarly. This construction can be iterated to prove the lemma.\looseness=-1
\end{proof}

It follows from Lemma \ref{linearlyindependent} that when $x\in \mathrm{Sys}(C)$ and $\mathrm{Sys}(C)$ is a locally top-dimensional stratum of $\mathcal{P}_{g}$ with $\{\nabla L(c)(x) \mid c\in C\}$ linearly independent, $C$ is a minimal filling set of curves. Also under these assumptions, by Corollary \ref{linear} and local finiteness the cone of increase at $x$ has nonempty intersection with the tangent cone to $\mathcal{P}_{g}$ at $x$.

A stratum will be said to be \textit{above} a point $x$ if the stratum contains a sequence $\{x_{i}\}$ of points converging to $x$ such that $f_{\mathrm{sys}}(x_{i})>f_{\mathrm{sys}}(x)$. Similarly, a level set of $f_{\mathrm{sys}}$ is above $x$ if the value of $f_{\mathrm{sys}}$ on this level set is larger than $f_{\mathrm{sys}}(x)$.

Recall that $\mathrm{Sys}(C)\subset E(C)$ for any set of systoles $C$.

\begin{Definition}[accidental systole]
Suppose there exists an open set $U$ in $\mathcal{T}_{g}$ such that $\mathrm{Sys}(C)\cap U$ is given by $E(C')\cap U$ for $C'\subsetneq C$ and $f_{\mathrm{sys}}$ is nonconstant on $\mathrm{Sys}(C)\cap U$. The curves in $C\setminus C'$ will be called accidental systoles at points of $\mathrm{Sys}(C)\cap U$.
\end{Definition}

The existence of accidental systoles means the linear independence condition in Lemma \ref{linearlyindependent} breaks down.

\begin{Remark}
\label{stillbalanced}
In the presence of accidental systoles, the vector $v_{C'}$ can be set equal to $v_{C}$, and since the convex hull of $\{\nabla L(c)(x) \mid c\in C'\}$ is contained in the convex hull of $\{\nabla L(c)(x) \mid c\in C\}$, if $E(C')$ is balanced at $x$, this implies that $\mathrm{Sys}(C)$ is balanced at $x$.
\end{Remark}

\begin{Remark}
The author is not aware of any examples of accidental systoles. Possibly they do not exist.
\end{Remark}

\begin{Theorem}
\label{balancedthm2}
Suppose a point $x\in \mathcal{P}_{g}$ has the property that the cone of increase of $f_{\mathrm{sys}}$ at $x$ is full. Then the intersection of this cone of increase with the tangent cone of $\mathcal{P}_{g}$ is nonempty.
\end{Theorem}
\begin{proof}
The theorem is proven by constructing a stratum above $x$ in $\mathcal{P}_{g}$. The first step is to show the existence of a nonzero vector $v_{C'}(x)$ in the cone of increase, where $C'$ is a filling subset of $C$, and for every $c\in C\setminus C'$, $L(c)$ is increasing faster in the direction of $v_{C'}(x)$ than the length of any curve $c'\in C'$. Suppose to begin with that $x$ is not a minimum of $C$ in $E(C)$.

Let $n$ be the dimension of $\mathrm{Span}\{\nabla L(c)(x) \mid c\in C\}$. Identify the subspace of $T_{x}\mathcal{T}_{g}$ spanned by~$\{\nabla L(c)(x) \mid c\in C\}$ with $\mathbb{R}^{n}$ via a vector space isomorphism that preserves the inner product on~$T_{x}\mathcal{T}_{g}$ induced by $\mathbf{g}$. Find a hyperplane in $\mathbb{R}^{n}$ containing the points represented by at least~$n$ of the vectors $\{\nabla L(c)(x) \mid c\in C\}$ and with the property that all other vectors in the set~$\{\nabla L(c)(x) \mid c\in C\}$ represent points on the opposite side of the hyperplane from the zero vector. The vector in $T_{x}\mathcal{T}_{g}$ corresponding to the closest point on the hyperplane to the zero vector is $v_{C'}(x)$. It follows from Lemma \ref{neededtofill} that $C'$ fills. The set of curves $C'$ cannot be eutactic at $x$, as this would contradict the assumption that the cone of increase at $x$ is full. Consequently, $v_{C'}(x)$ is not zero.

When $\{\nabla L(c)(x) \mid c\in C'\}$ are not linearly independent, as a result of second and higher order rates of change of lengths, it does not necessarily follow that $v_{C'}(x)$ is in the tangent cone to a~stratum $\mathrm{Sys}(C')$ at $x$. In order to show that there does, nevertheless, exist a stratum in~$\mathcal{P}_{g}$ above~$x$, a ``Voronoi decomposition'' of the unit tangent cone $UT_{x}\mathcal{T}_{g}$ will be defined, and properties of this decomposition invariant under arbitrarily small deformations in length will be used.

The decomposition\footnote{It will not be shown that the decompositions used in this proof truly are cell decompositions, as this is cumbersome, and not essential. If the reader is concerned, it is possible to simply assume that unspecified subdivisions have been performed to obtain cell decompositions.} $\mathcal{V}(C,x)$ of $UT_{x}\mathcal{T}_{g}$ has a top-dimensional cell labelled by a curve $c\in C$ consisting of all $v\in UT_{x}\mathcal{T}_{g}$ for which
\begin{equation*}
\mathbf{g}(v, -\nabla L(c)(x))>\mathbf{g}(v, -\nabla L(c')(x)) \, \forall c'\in C\setminus\{c\}.
\end{equation*}
Smaller-dimensional cells of the decomposition are defined analogously, and labelled by subsets of~$C$. Note that it follows from Proposition~\ref{balancedlemnonfilling} and Remark \ref{subsetrem} that every curve in $C$ determines a top-dimensional cell of $\mathcal{V}(C,x)$.

In what follows, it will only be necessary to study the intersection of the cell decomposition with the unit vectors in the span of $\{\nabla L(c)(x) \mid c\in C\}$. This will be denoted by $\mathcal{V}(C,x)^{\mathrm{Span}}$. It follows from Lemma~\ref{neededtofill} that zero-dimensional cells of $\mathcal{V}(C,x)^{\mathrm{Span}}$ are labelled by filling subsets of $C$. It was shown that the intersection of $\mathcal{V}(C,x)^{\mathrm{Span}}$ with the cone of increase contains a~0-dimensional cell labelled by $C'$.

When the dimension of $\mathrm{Span}\{\nabla L(c)(x) \mid c\in C\}$ stays constant on a neighbourhood of $x$, define~$\mathcal{V}(C,x,\epsilon)^{\mathrm{Span}}$ as follows: scale the unit vectors in $UT_{x}\mathcal{T}_{g}\cap \mathrm{Span}\{\nabla L(c)(x) \mid c\in C\}$ by $\epsilon$, and take the image under the exponential map to obtain an embedded $n$-dimensional sphere $\mathcal{S}_{\epsilon}$. Assume~$\epsilon$ is small enough that local finiteness ensures that the systoles everywhere on $\mathcal{S}_{\epsilon}$ are contained in the set $C$. $\mathcal{V}(C,x,\epsilon)^{\mathrm{Span}}$ is the cell decomposition of $\mathcal{S}_{\epsilon}$ induced by the stratification; a cell is labelled by a subset $C''$ of $C$ if the points of the cell are contained in $\mathrm{Sys}(C'')$.

\begin{figure}[ht]
\centering
\includegraphics[width=0.6\textwidth]{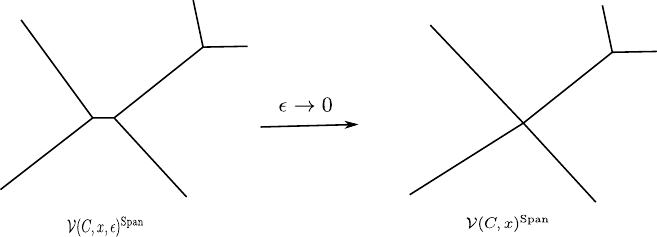}
\caption{A cell of small diameter in $\mathcal{V}(C,x,\epsilon)^{\mathrm{Span}}$ can shrink to a point as $\epsilon$ is decreased to zero.}
\label{Voronoifig}
\end{figure}

There is a sense in which $\mathcal{V}(C,x,\epsilon)^{\mathrm{Span}}$ can be made arbitrarily close to $\mathcal{V}(C,x)^{\mathrm{Span}}$ for sufficiently small $\epsilon$. For small $\epsilon$, the cells of $\mathcal{V}(C,x,\epsilon)^{\mathrm{Span}}$ and $\mathcal{V}(C,x)^{\mathrm{Span}}$ correspond except when a cell of small diameter in $\mathcal{V}(C,x,\epsilon)^{\mathrm{Span}}$ shrinks to a point as $\epsilon$ is decreased to zero. A~lower-dimensional analogue of this is illustrated in Figure~\ref{Voronoifig}. In particular, there exists a~0-dimensional cell of $\mathcal{V}(C,x,\epsilon)^{\mathrm{Span}}$ corresponding to a vector in the cone of increase of $f_{\mathrm{sys}}$ at $x$. This shows that at $x$, there is nonempty intersection of the cone of increase with the tangent cone of $\mathcal{P}_{g}$, as required.

When the dimension of $\mathrm{Span}\{\nabla L(c) \mid c\in C\}$ does not stay constant on a neighbourhood of~$x$, i.e., the dimension of the span increases, defining $\mathcal{V}(C,x,\epsilon)^{\mathrm{Span}}$ involves constructing a~cross section with tangent space that will be identified with $\mathrm{Span}\{\nabla L(c)(x) \mid c\in C\}$. A cross section will be used that does not ``see'' the increase in the dimension of the span. The purpose of this is to obtain a stratum labelled by a subset of $C$ the elements of which have lengths whose gradients at $x$ have the same span as $\{\nabla L(c)(x) \mid c \in C\}$. On a neighbourhood of $x$ in $\mathcal{T}_{g}$, choose $\mathrm{Croppedspan}\{\nabla L(c) \mid c\in C\}$ to satisfy the following conditions:
\begin{itemize}\itemsep=0pt
\item For some fixed $C'\subset C$ independent of $y$, $\mathrm{Croppedspan}\{\nabla L(c)(y) \mid c\in\! C\}$ is an $n$-dimensional subspace of $T_{y}\mathcal{T}_{g}$ spanned by $\{\nabla L(c)(y) \mid c\in C'\}$.
\item{When $y=x$, $\mathrm{Croppedspan}\{\nabla L(c)(y) \mid c\in C\}=\mathrm{Span}\{\nabla L(c)(x) \mid c\in C\}$.}
\item The projection of $\nabla L(c)(y)$ to the orthogonal complement of $\mathrm{Croppedspan}\{\nabla L(c)(y) \mid c\in C\}$ in $T_{y}\mathcal{T}_{g}$ for every $c\in C$ has length bounded from above by a function of distance from~$x$ that approaches zero as $y$ approaches~$x$.
\end{itemize}
The existence of $\mathrm{Croppedspan}\{\nabla L(c) \mid c\in C\}$ will now be shown. Choose and fix any $C'\subset C$ satisfying the first two conditions. The third condition then follows from the fact that every curve $c\in C\setminus C'$ has length whose gradient becomes linearly dependent on $\{\nabla L(c)(x) \mid c\in C\}$ at~$x$.

Suppose $\mathcal{B}'_{\alpha}$ has the property that $\mathcal{B}'_{\alpha}\setminus x$ is foliated by $\mathcal{S}_{\epsilon}$ for $\epsilon\in (0,\alpha]$, and let $\phi\colon \mathcal{D}^{n}\rightarrow \mathcal{B}'_{\alpha}$ be a diffeomorphism of the $n$-dimensional unit ball $\mathcal{D}^{n}$ into $\mathcal{B}'_{\alpha}$. There is an isotopy of $\phi$ to a map with image a ball $\mathcal{B}_{\alpha}$ centered on $x$ with tangent space at every point $y$ given by $\mathrm{Croppedspan}\{\nabla L(c)(y) \mid c\in C\}$. Then $\mathcal{V}(C,x,\epsilon)^{\mathrm{Span}}$ is defined to be the cell decomposition induced by the stratification on the sphere of radius $\epsilon$ in $\mathcal{B}_{\alpha}$.

Using $\mathcal{V}(C,x,\epsilon)^{\mathrm{Span}}$, it is then possible to argue as in the previous case that at $x$, there is nonempty intersection of the cone of increase with the tangent cone of $\mathcal{P}_{g}$. This is because there is a 0-dimensional cell in $\mathcal{V}(C,x,\epsilon)^{\mathrm{Span}}$ above $x$ labelled by a subset $C'$ of $C$ for which~${\mathrm{Span}\{\nabla L(c)(x) \mid c\in C'\}=\mathrm{Span}\{\nabla L(c)(x) \mid c\in C\}}$.

It remains to prove the theorem in the case that $x$ is a minimum of $C$ on $E(C)$ and the cone of increase at $x$ is full. Let $x_{1}, x_{2}, \dots$ be a sequence of points in $\mathcal{T}_{g}$ with $\lim_{i\rightarrow \infty}x_{i}=x$ and such that no $x_{i}$ is a minimum of $C$. The sequence is also assumed to be contained within a neighbourhood $\mathcal{N}$ of $x$ on which the systoles are all contained in $C$ and every point has an open cone of directions in which the lengths of all curves in $C$ are strictly increasing. Within~$\mathcal{N}$, define a relative stratum $E^{\mathrm{rel}}(C', d)$ to be the set of all points on $E(C', d)$ at which
\begin{equation*}
L(c')+d|_{c'}<L(c)+d|_{c} \forall c'\in C'\qquad\text{and}\qquad\forall c\in C\setminus C'.
\end{equation*}
The set of curves $C'$ will be called relative systoles. The previous arguments show the existence of a relative stratum $E^{\mathrm{rel}}(C'(x_{i}), d(x_{i}))$ above each point $x_{i}$, where $C'(x_{i})\subset C$ is a filling set, and $d(x_{i})$ is the value of $d$ for which $E(C, d(x_{i}))$ contains $x_{i}$. Assume, in addition, that the points of the sequence have been chosen such that they all have relative strata given by the same subsets of $C$. Every point on such a relative stratum in the interior of $\mathcal{N}$ has another relative stratum with filling relative systoles above it, defined using the tuple $d(x_{i})$. In the limit as $x_{i}$ approaches $x$, $d(x_{i})$ approaches $d(x)=0$, and these relative strata limit to strata in $\mathcal{P}_{g}$ above $x$. This concludes the proof of the theorem.
\end{proof}

\begin{Theorem}
\label{balancedthm}
Where a stratum of $\mathcal{P}_{g}$ is locally top-dimensional it is balanced. If $x\in \mathcal{P}_{g}$ is not balanced, there are points directly above $x$ in $\mathcal{P}_{g}$ that are balanced.
\end{Theorem}

\begin{proof}
The theorem is proven by assuming $x$ is in a locally top-dimensional stratum $\mathrm{Sys}(C)$ of~$\mathcal{P}_{g}$ and using Proposition \ref{plus1} to derive a contradiction when $\mathrm{Sys}(C)$ is not balanced at $x$.

Suppose $x\in \mathrm{Sys}(C)$ is in a locally top-dimensional stratum in $\mathcal{P}_{g}$ that is not balanced at $x$. It follows that $x$ is not a critical point of $f_{\mathrm{sys}}$. If $x$ were a boundary point of $f_{\mathrm{sys}}$, as discussed in Section \ref{criticalstructure}, these are known to have higher-dimensional balanced strata in $\mathcal{P}_{g}$ above them. It is therefore possible to assume that the cone of increase at $x$ is full.

On a neighbourhood of $x$, $\mathrm{Sys}(C)$ is not contained in a level set of $f_{\mathrm{sys}}$. If it were, Theorem~\ref{balancedthm2} would give a stratum $\mathrm{Sys}(C')$ above $x$ in $\mathcal{P}_{g}$ with $C'\subset C$ and for which there exists a~$v_{C'}$ in the tangent cone to $\mathrm{Sys}(C')$ at points sufficiently close to $x$. This shows it is possible to construct a stratum in $\mathcal{P}_{g}$ directly above $x$ and with larger dimension than $\mathrm{Sys}(C)$, contradicting the assumption that $x$ is in a locally top-dimensional stratum.

Denote by $c(f)$ the set of curves corresponding to the vertices of $f$, where $f$ is the face of the convex hull of $\{\nabla L(c)(x) \mid c\in C\}$ from Proposition \ref{plus1}. The existence of a larger-dimensional stratum above $x\in \mathrm{Sys}(C)$ is a consequence of Proposition \ref{plus1} and the fact that every facet~$f_{i}$ of the convex hull of $\{\nabla L(c)(x) \mid c\in C\}$ determines a stratum $\mathrm{Sys}(c'(f_{i}))$ adjacent to $x$, where~${c'(f_{i})\subset c(f_{i})}$ fills the same subsurface of $\mathcal{S}_{g}$ as $c(f_{i})$. To prove this fact about facets, note that the span of $\{\nabla L(c)(x) \mid c\in c(f_{i})\}$ has a 1-dimensional orthogonal complement in~$\mathrm{Span}\{\nabla L(c)(x) \mid c\in C\}$. Denote by $o_{f_{i}}$ a tangent vector in $T_{x}\mathcal{T}_{g}$ to this orthogonal complement. The fact that $f_{i}$ is a facet is used to ensure that there is a choice of~$o_{f_{i}}$ with strictly positive inner product with every element of~$\{\nabla L(c)(x) \mid c\in C\setminus c(f_{i})\}$. To first order, $o(f_{i})$ appears to be in the tangent cone at $x$ to a stratum $\mathrm{Sys}(c(f_{i}))$ adjacent to $x$. The construction in the proof of Theorem \ref{balancedthm2} then shows that $o_{f_{i}}$ is in the tangent cone at $x$ to a stratum~$\mathrm{Sys}(c'(f_{i}))$, where~$c'(f_{i})\subset c(f_{i})$ fills the same subsurface as $c(f_{i})$. The stratum $\mathrm{Sys}(c'(f_{i}))$ has higher dimension than $\mathrm{Sys}(C)$ at $x$, giving the required contradiction.
\end{proof}

\begin{Remark}
\label{thedifficulty}
There are two obstacles in trying to prove that sets of systoles on locally top-dimensional strata of $\mathcal{P}_{g}$ are all minimal filling sets. The first one is in ruling out accidental systoles. The second is that the converse of Lemma \ref{lemma14} does not hold. Not every filling subset $C'$ of $C$ determines a set of minima $\mathrm{Min}(C')$ contained in $\partial \mathrm{Min}(C)$. This could give rise to higher-dimensional analogues of the following: Let $C=\{c_{1}, c_{2}, c_{3}, c_{4}\}$ where $\{\nabla L(c)(x) \mid c\in C\}$ are linearly dependent and $v_{C}(x)$ is nonzero and in the convex hull of $\{\nabla L(c)(x) \mid c\in C\}$. Suppose also that $\{c_{1}, c_{2}\}$ is a minimal filling set, and $c_{1}$ and $c_{2}$ are the labels of a pair of vertices of the convex hull of $\{\nabla L(c)(x) \mid c\in C\}$ diagonally opposite each other. There seems to be no reason to expect that any of the faces of the convex hull determine a set of filling systoles.
\end{Remark}

\section{Critical points and boundary points of the systole function}
\label{criticalstructure}
The purpose of this section is to study the structure of $\mathcal{P}_{g}$ around critical and boundary points of $f_{\mathrm{sys}}$. Section \ref{subexample} later gives an example that is interesting in its own right, illustrating the theorem proven in this section. The exposition begins with critical points.

Recalling the definition of a topological Morse function, on a neighbourhood of a critical point~$p$ of index $j$ there is a homeomorphism $\psi$ that fixes the critical point. This homeomorphism has the property that there exists a smooth chart on a neighbourhood of $p$ such that
\[
f_{\mathrm{sys}}\circ\psi(x)-f_{\mathrm{sys}}(p)=-x_{1}^{2}-\dots-x_{j}^{2}+x_{j+1}^{2}+\dots+x_{n}^{2}.
\]
In this sense, the level sets of $f_{\mathrm{sys}}$ on a neighbourhood of a critical point of index $j$ are the same, up to homeomorphism, as for a critical point of index $j$ of a smooth Morse function. It was shown in \cite{Akrout} that critical points of $f_{\mathrm{sys}}$ occur exactly where $\mathrm{Sys}(C)$ intersects $\mathrm{Min}(C)$. The index $j$ is then equal to the dimension of the span of $\{\nabla L(c) \mid c\in C\}$ at the critical point.

\begin{figure}[h!]
\centering
\includegraphics[width=0.5\textwidth]{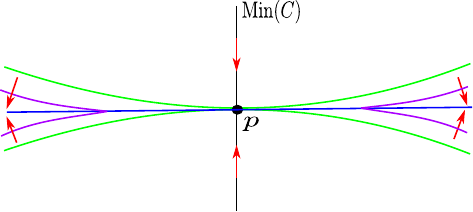}
\caption{The level sets of the systole function passing through the critical point are shown in green. A level set at higher value of $f_{\mathrm{sys}}$ is shown in violet. The red arrows show the direction in which $f_{\mathrm{sys}}$ is increasing. The simplices $\sigma_{1}$ and $\sigma_{2}$ are shown in blue.}
\label{localstructurefig}
\end{figure}

Suppose $p$ is a critical point of index $j$ in the intersection of $\mathrm{Sys}(C)$ with $\mathrm{Min}(C)$. It follows from local finiteness that there is a neighbourhood of $p$ in $\mathcal{T}_{g}$ on which the systoles are contained in the set $C$. Convexity of length functions implies that the subspace $\{\nabla L(c) \mid c\in C\}^{\perp}$ is contained in the cone of increase of $f_{\mathrm{sys}}$ at $p$. Within $\mathrm{Min}(C)$, in every direction there is by definition a curve or curves in $C$ whose length is decreasing away from $p$ to first order. For a~point~$q$ in $\mathrm{Min}(C)$ near $p$, the systoles are the subset of curves in $C$ whose lengths are decreasing fastest away from $p$ in the direction of $q$. It follows that $j\leq |C|-1$ and on a neighbourhood of $p$, the level set of $f_{\mathrm{sys}}$ containing $p$ intersects $\mathrm{Min}(C)$ in the single point $p$. This is illustrated schematically in Figure~\ref{localstructurefig} for the toy example with two systoles.

For $p$ in $\overline{\mathrm{Sys}(C)}$, the tangent cone to $\mathrm{Sys}(C)$ at $p$ will be denoted by $TC_{p}\mathrm{Sys}(C)$.

\begin{Proposition}
\label{localstructure}
Fix a critical point $p\in \mathcal{P}_{g}$ of index $j$ with set of systoles $C$. A neighbourhood of $p$ in $\mathcal{P}_{g}$ contains a piecewise smooth cell of codimension $j$ on which $f_{\mathrm{sys}}$ is increasing to second order away from the critical point. The tangent space of this cell at $p$ is $\{\nabla L(c)(p) \mid c\in C\}^{\perp}$. The subspace spanned by $\{\nabla L(c)(p) \mid c\in C\}$ consists of directions in which the systole function is decreasing away from $p$.
\end{Proposition}

\begin{proof}
By local finiteness, convexity of length functions and the definition of $\mathrm{Min}(C)$, the subspace spanned by $\{\nabla L(c)(p) \mid c\in C\}$ consists of directions in which $f_{\mathrm{sys}}$ in decreasing away from $p$ to first order and all the directions in which $f_{\mathrm{sys}}$ is increasing at $p$ to second order are exactly given by $\{\nabla L(c)(p) \mid c\in C\}^{\perp}$.

The subspace $\{\nabla L(c)(p) \mid c\in C\}^{\perp}$ of $T_{p}\mathcal{T}_{g}$ is tangent to the level set of $f_{\mathrm{sys}}$ containing $p$. Recall that there exists a triangulation of $\mathcal{P}_{g}$ compatible with the stratification. Denote by $\sigma=\bigcup_{i=1}^{n}\sigma_{i}$ the union of simplices for which $TC_{p}\sigma=\{\nabla L(c)(p) \mid c\in C\}^{\perp}$. The existence of $\sigma$ is guaranteed by two observations: firstly, every vector in the cone of increase at $p$ is in $\{\nabla L(c) \mid c\in C\}^{\perp}$ and hence tangent to the level sets passing through $p$. Secondly, the restriction of $f_{\mathrm{sys}}$ to the interior of a stratum or a simplex is smooth, so the restriction of the level sets to the interiors of the simplices can only contain the smooth pieces of the piecewise smooth level sets. This is illustrated schematically in Figure~\ref{localstructurefig}.

It will now be shown that, on a neighbourhood of $p$, $\sigma$ is contained in $\mathcal{P}_{g}$. This starts by showing that each of the locally top-dimensional simplices $\sigma_{1}, \dots, \sigma_{n}$ of $\sigma$ is necessarily contained in locus that is balanced near $p$. Let $x$ be a point in the interior of a simplex $\sigma_{i}$ in $\sigma$ where the systoles on $\sigma_{i}$ are given by $C_{i}\subset C$, and let $V$ be the subspace of $T_{x}\mathcal{T}_{g}$ given by the orthogonal complement of $v_{C_{i}}(x)$. From the known geometry of the level sets near $p$, it follows that any vector in $V$ is either tangent to a level set of $f_{\mathrm{sys}}$ restricted to $\sigma_{i}$ or to a geodesic $\gamma$ with $\gamma(0)=x$ passing through level sets below $x$ for $t\in (-\epsilon,\epsilon)$. In Figure~\ref{localstructurefig}, $\gamma$ is approximately parallel to the~$y$ axis. Local finiteness then implies that $C_{i}$ is $V$-eutactic, and hence $\sigma_{i}$ is balanced.

Consider a sequence of points in $\sigma_{i}$ approaching $p$ from above. In the limit as $p$ is approached, $f_{\mathrm{sys}}|_{\sigma_{i}}$ is increasing only to second order. At points in the sequence, the balanced sets of gradients~$\{\nabla L(c) \mid c\in C_{i}\}$ therefore become arbitrarily close to being eutactic. This implies that~$C_{i}$ is a filling set, since otherwise a contradiction is obtained by invoking Proposition \ref{Thurstonprop} and observing that in the limit as $p$ is approached, the open cone of directions in which the lengths of curves in $C_{i}$ are increasing is empty. Since $\mathcal{P}_{g}$ is closed, this concludes the proof that $\sigma\subset \mathcal{P}_{g}$.
\end{proof}

\begin{Remark}
Proposition \ref{localstructure} implies that for a critical point $p\in \mathrm{Sys}(C)$,
codimension of $\mathcal{P}_{g} \leq |C|-1$ or dimension of $\mathcal{P}_{g}\geq 6g-|C|-5$.
This does not require Thurston's construction from \cite{Thurston} but as explained above, is a consequence of Akrout's characterisation of critical points of $f_{\mathrm{sys}}$ and Proposition \ref{Thurstonprop}. This has not been well understood in the literature, with conjectures made and lengthy proofs given of special cases of this statement. When $p$ is not a critical point, but is contained in a locally top-dimensional stratum of $\mathcal{P}_{g}$, by Theorem \ref{balancedthm}, the cardinality of the systole set similarly gives an upper bound on the codimension of the stratum and hence of~$\mathcal{P}_{g}$.\looseness=-1
\end{Remark}

\textbf{Boundary points of the systole function and their structure.} Recall that a boundary point of $f_{\mathrm{sys}}$ is a point at which $\mathrm{Sys}(C)$ intersects $\mathrm{Min}(C')$, for $C'\subsetneq C$. Since $\mathrm{Min}(C')$ is empty unless $C'$ fills, boundary points of $f_{\mathrm{sys}}$ are all necessarily contained in $\mathcal{P}_{g}$. Boundary points of~$f_{\mathrm{sys}}$ were also shown in \cite{SchmutzMorse} to be isolated and hence there are at most finitely many modulo the action of $\Gamma_{g}$.

\begin{figure}
\centering
\includegraphics[width=0.5\textwidth]{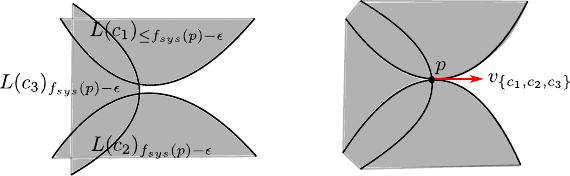}
\caption{The intersection of level sets near a boundary point $p$ of $f_{\mathrm{sys}}$. The vector $v_{\{c_{1}, c_{2}, c_{3}}$ shows a~direction in which $f_{\mathrm{sys}}$ is increasing away from $p$ to second order. }
\label{sliver}
\end{figure}

A boundary point of $f_{\mathrm{sys}}$ is not a critical point, but it also has no open cone in which $f_{\mathrm{sys}}$ is increasing. The level sets of the curves in $C$ intersect as illustrated in Figure~\ref{sliver}. The function~$f_{\mathrm{sys}}$ increases away from a boundary point of $f_{\mathrm{sys}}$ only to second order.

Ignoring the curves in $C\setminus C'$, around $p$, it is possible to construct the same simplices $\sigma_{1}, \dots, \sigma_{n}$ as in Proposition \ref{localstructure}. For any vector in the tangent cone to $\sigma$ at $p$ giving a direction in which the length of a curve in $C\setminus C'$ is decreasing away from $p$, the set of systoles is in $C\setminus C'$. In all other directions near $p$, the behaviour of $f_{\mathrm{sys}}$ is identical to that near a critical point. Note that all directions in which $f_{\mathrm{sys}}$ is increasing away from $p$ are in the tangent cone to $\mathcal{P}_{g}$ at $p$.
\section{Thurston's vector field}
\label{Tvf}
In the 80s, Thurston gave a talk at Princeton on the construction of a $\Gamma_{g}$-equivariant deformation retraction of $\mathcal{T}_{g}$ onto $\mathcal{P}_{g}$. This was done by showing the existence of a vector field that generates a flow used in the construction of the deformation retraction. Notes from this talk are given in~\cite{Thurston}. The purpose of this section is to explain what Theorems \ref{balancedthm}, \ref{balancedthm2} and Proposition \ref{localstructure} imply about this vector field.

The vector field from \cite{Thurston} will first be defined. Choose and fix an $\epsilon>0$. For \textit{finite} $C$, define
\begin{equation*}
U_{C}=\{x\in \mathcal{T}_{g} \mid C \text{ is exactly the set of curves of length less than }f_{\mathrm{sys}}+|C|\epsilon\},
\end{equation*}
where $|C|$ is the cardinality of the set $C$. It was shown that when $\epsilon$ is sufficiently small, the set
\begin{equation}
\label{cover}
\{U_{C} \mid C \text{ is a finite set}\}
\end{equation}
covers $\mathcal{T}_{g}$. This uses the well-known result (see, for example, \cite{M}) that the number $n(L)$ of curves on a hyperbolic surface corresponding to a point $x\in \mathcal{T}_{g}$ of length less than or equal to $L$ increases faster than linearly with $L$. Consequently, $\epsilon$ must be made small otherwise a set $C$ of curves of length less than $f_{\mathrm{sys}}(x)+|C|\epsilon$ quickly becomes infinite, and the set in equation~\eqref{cover} is not a~cover.

On $U_{C}$, for nonfilling $C$, a vector field $X_{C}$ is defined to be the vector of unit length such that
\begin{equation*}
\sum_{c\in C}\log X_{C}(L(c))
\end{equation*}
is maximised. When $C$ fills, $X_{C}$ is defined to be zero. Note that this formula assumes real logs, i.e., $X_{C}$ is a derivation whose evaluation on each of the functions $\{L(c) \mid c\in C\}$ is strictly positive. A smooth vector field $X_{\epsilon}$ was obtained with the help of a partition of unity subordinate to the cover. This gives a vector field, nonvanishing away from a neighbourhood of $\mathcal{P}_{g}$, and in a direction of increasing $f_{\mathrm{sys}}$. For any $x\in \mathcal{T}_{g}\setminus \mathcal{P}_{g}$, by choosing $\epsilon$ sufficiently small, it can be ensured that $X_{\epsilon}(x)\neq 0$. More details on the construction of $X_{\epsilon}$ can be found in \cite{Thurston}.

\begin{Definition}[gap function $g_{\mathrm{sys}}$]\label{Def8}
The \textit{gap function} $g_{\mathrm{sys}}\colon\mathcal{T}_{g}\rightarrow \mathbb{R}_{+}$ is given by
\begin{equation*}
x\mapsto \inf\{r\in \mathbb{R} \mid \text{the curves of length at most }f_{\mathrm{sys}}(x)+r\text{ fill}\}.
\end{equation*}
\end{Definition}

By local finiteness, $g_{\mathrm{sys}}$ is zero only on $\mathcal{P}_{g}$, and increases away from $\mathcal{P}_{g}$ on some neighbourhood of $\mathcal{P}_{g}$.

\begin{Remark}
Close to $\mathcal{P}_{g}$, $g_{\mathrm{sys}}$ behaves like a measure of distance from $\mathcal{P}_{g}$. Globally it is not an easy function to work with. For example, $g_{\mathrm{sys}}$ is not a topological Morse function. This follows from the observation that at some points, the cone of increase of $g_{\mathrm{sys}}$ is obtained by taking a~\textit{union} of cones of increase corresponding to difference choices of longest curve in a filling set. The cone of increase can therefore become disconnected when the zero vector is removed.
\end{Remark}

\begin{Proposition}
\label{atworstparallel}
There is a neighbourhood $\mathcal{N}_{g}$ of $\mathcal{T}_{g}$ such that for every point $x\in \mathcal{N}_{g}\setminus \mathcal{P}_{g}$, choosing $\epsilon$ small enough to ensure that $X_{\epsilon}(x)\neq 0$ ensures that $g_{\mathrm{sys}}$ is strictly decreasing in the direction of $X_{\epsilon}(x)\neq 0$.
\end{Proposition}
\begin{proof}
The neighbourhood $\mathcal{N}_{g}$ of $\mathcal{P}_{g}$ is required to satisfy the following conditions
\begin{itemize}\itemsep=0pt
\item $\mathcal{N}_{g}$ is a regular, $\Gamma_{g}$-equivariant neighbourhood of $\mathcal{P}_{g}$.
\item $\mathcal{N}_{g}$ is small enough so that $g_{\mathrm{sys}}$ is increasing radially away from $\mathcal{P}_{g}$ on $\mathcal{N}_{g}$.
\item $\mathcal{N}_{g}$ is small enough so that local finiteness guarantees that the systoles on $\mathcal{P}_{g}$ determine the systoles within $\mathcal{N}_{g}$.
\end{itemize}

For sufficiently small $\mathcal{N}_{g}$, near critical points and boundary points of $f_{\mathrm{sys}}$, the proposition follows from Proposition \ref{localstructure} and the characterisation of boundary points of $f_{\mathrm{sys}}$ given at the end of Section \ref{criticalstructure}. Away from critical points and boundary points, for sufficiently small $\mathcal{N}_{g}$ the proposition follows from Theorems \ref{balancedthm2} and \ref{balancedthm}. In $\mathcal{N}_{g}$ the gradients of the systoles in the top-dimensional strata of $\mathcal{T}_{g}$ point radially inwards towards $\mathcal{P}_{g}$, and the same is then true of the linear combinations of these gradients where the top-dimensional strata meet along lower-dimensional strata. Where cells of $\mathcal{P}_{g}$ come together as shown schematically in Figure~\ref{wedgie}, Theorem \ref{balancedthm} rules out examples such as in the left hand side of Figure~\ref{wedgie}, where there is a conflict between increasing $f_{\mathrm{sys}}$ and decreasing $g_{\mathrm{sys}}$.
\begin{figure}
\centering
\includegraphics[width=0.5\textwidth]{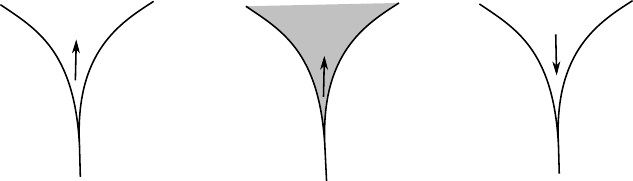}
\caption{The black lines/shaded region represent points on $\mathcal{P}_{g}$, and the arrows represent the direction of increasing $f_{\mathrm{sys}}$. Where cells of $\mathcal{P}_{g}$ come together, Theorem \ref{balancedthm} ensures there is no conflict between increasing $f_{\mathrm{sys}}$ and decreasing $g_{\mathrm{sys}}$; the behaviour schematically represented by the figure on the left is not possible.}
\label{wedgie}
\end{figure}
Theorems \ref{balancedthm2} and \ref{balancedthm} also rule out a vector field pointing away from the boundary of $\mathcal{P}_{g}$ at $x$.
\end{proof}

\section{Example: Schmutz's critical point}
\label{subexample}
The purpose of this subsection is to discuss an application of Proposition \ref{localstructure} in understanding the geometry of $\mathcal{P}_{2}$. This critical point is due to Schmutz \cite{SchmutzMorse} and is the first of a family of examples, one in each genus, of critical points of coindex equal to the virtual cohomological dimension of the mapping class group. This family of examples plays a key role in the study of the Steinberg module, \cite{Steinberg}.

If $p$ is a critical point of index $j$, and $C$ is the set of systoles at $p$, a \textit{descendent} $C'$ of $C$ is a~subset of $C$ with the following properties:
\begin{enumerate}\itemsep=0pt
\item[(1)] $C'$ is eutactic at $p$.
\item[(2)] The span of $\{\nabla L(c) \mid c\in C'\}$ is equal to the span of $\{\nabla L(c) \mid c\in C\}$.
\item[(3)] $C'$ does not have a proper subset satisfying the first two conditions.
\end{enumerate}
Note that the cardinality of a descendent $C'$ is $j+1$.

\begin{figure}
\centering
\includegraphics[width=0.7\textwidth]{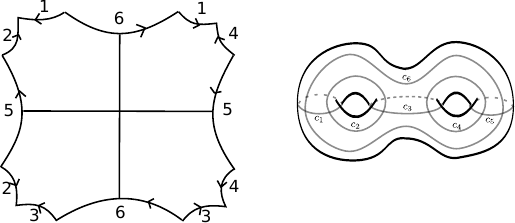}
\caption{The set of 6 systoles at the critical point in the example are shown on the right. The left side of the figure shows a fundamental domain with edges lying along the systoles, and with numbers indicating the glung maps (not the labels of the curves). This figure is taken from \cite{Steinberg}.}
\label{genus2example1}
\end{figure}

Let $C$ be the set of curves on the right hand side of Figure~\ref{genus2example1}. It was shown in \cite[Theorem~36]{SchmutzMorse} that there is a critical point $p$ of index 3 at which this set of curves is realised as the set of systoles. At $p$, the systoles intersect at right angles, cutting the surface into four regular hexagons. In~\cite{SchmutzMorse}, it was shown that $\mathrm{Min}(C)$ is a cell of dimension 3 with empty boundary.

Deleting any pair of intersecting curves in $C$ gives a minimal filling set of curves, all of which are descendents of $C$. The complete set of descendents is
\begin{gather*}
C_{1}=\{c_{1}, c_{2}, c_{3}, c_{4}\}, \qquad C_{2}=\{c_{2}, c_{3}, c_{4}, c_{5}\}, \qquad C_{3}=\{c_{3}, c_{4}, c_{5}, c_{6}\},\\
C_{4}=\{c_{4}, c_{5}, c_{6}, c_{1}\}, \qquad C_{5}=\{c_{5}, c_{6}, c_{1}, c_{2}\}, \qquad C_{6}=\{c_{6}, c_{1}, c_{2}, c_{3}\}.
\end{gather*}

For each descendent $C_{i}$, $i=1, \dots, 6$, the locus $E(C_{i})$ is an embedded submanifold of dimension~3, intersecting $\mathrm{Min}(C)=\mathrm{Min}(C_{i})$ in the single point $p$, with $T_{p}E(C_{i})$ given by $\mathrm{Span}\{\nabla L(c) |\allowbreak c\in C_{i}\}^{\perp}=\mathrm{Span}\{\nabla L(c) \mid c\in C\}^{\perp}$. This follows from the results in Section \ref{Minimalfilling}, or on a neighbourhood of $p$ by the rank calculations in~\cite{SchmutzMorse}.

\begin{figure}
\centering
\includegraphics[width=0.6\textwidth]{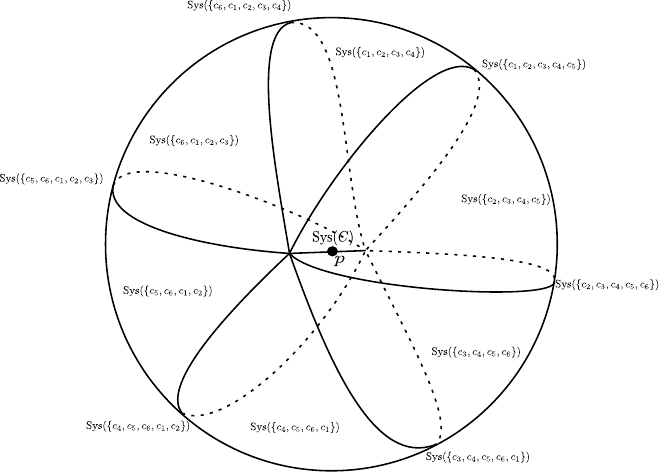}
\caption{The arrangement of strata in $\mathcal{P}_{2}$ on a neighbourhood of the critical point $p$ in the example. The 3-sphere is piecewise smooth, and has a symmetry group given by the dihedral group $D_{6}$ that acts transitively on the top-dimensional strata of $\mathcal{P}_{2}$ adjacent to $p$. }
\label{sphereexample}
\end{figure}

Near $p$, $\mathcal{P}_{2}$ could potentially lie along any $E(C_{i})$; this is determined by second order behaviour. It is a consequence of the main construction in \cite{estimating} that $E(C)=\bigcap_{i=1}^{i=6}E(C_{i})$ has a~connected component consisting of a one-dimensional embedded submanifold of $\mathcal{T}_{2}$ containing~$p$. On a~neighbourhood of $p$, by local finiteness, the points of $E(C)$ are in $\mathrm{Sys}(C)$. The same construction from \cite{estimating} also implies that along $E(C)\setminus \{p\}$ the set of gradients $\{\nabla L(c) \mid c\in C\}$ is linearly independent. It follows from the argument in Lemma \ref{linearlyindependent} that every proper filling subset of $C$ determines a stratum of $\mathcal{P}_{2}$ adjacent to $E(C)$; these are arranged as shown in Figure~\ref{sphereexample}.

The automorphism group of the hyperbolic surface corresponding to $p$ contains a subgroup of the mapping class group isomorphic to the dihedral group $D_{6}$, and the locus $E(C)$ passing through $p$ is the fixed point set of this subgroup. It follows that this connected component of~$E(C)$ is a geodesic with respect to any mapping class group-equivariant metric on $\mathcal{T}_{2}$, such as the Teichm\"uller metric or the Weil--Petersson metric.

The set of minima $\mathrm{Min}(C)$ also has a simple parameterisation. In genus 2, $\mathrm{Min}(C)$ is the set of points at which the curves in $C$ intersect at right angles.

\begin{Theorem}[{\cite[Theorem 3.5.14]{Ratcliffe}}]
\label{Rat}
For any $a,b,c \in \mathbb{R}_{+}$, there is a right-angled hyperbolic convex hexagon, unique up to congruence, with alternate sides of length $a,b$ and $c$.
\end{Theorem}

Using Theorem \ref{Rat}, from the left side of Figure~\ref{genus2example1}, it is not difficult to verify that by Theorem~\ref{Rat} the lengths of 3 alternate edges of any one of the hexagons of $S_{g}\setminus C$ completely parameterise $\mathrm{Min}(C)$.

\subsection*{Acknowledgements}
The author would like to thank Stavros Garoufalidis, Don Zagier and a number of very helpful anonymous referees for comments that led to significant improvements.

\pdfbookmark[1]{References}{ref}
\LastPageEnding

\end{document}